\documentclass[12pt,a4paper]{article}

\usepackage[english]{babel}
\usepackage{verbatim}
\usepackage{bm}

\usepackage{amsmath}
\usepackage{amssymb}
\usepackage{amsfonts}
\usepackage{amsbsy}

\usepackage{graphicx}
\usepackage{graphics}
\usepackage{subfig}
\usepackage{cite}
\usepackage{array}
\usepackage{booktabs} 
\usepackage{caption}  
\usepackage{tabularx}

\usepackage{hyphenat} 

\usepackage{amsthm} 

\newtheorem*{proof*}{Proof}
\newtheorem*{remark*}{Remark}

\numberwithin{equation}{section}
\numberwithin{theorem}{section}
\numberwithin{lemma}{section}
\numberwithin{defin}{section}
\numberwithin{cor}{section}
\numberwithin{prop}{section}

\newcommand{\Q}{\mathbf{Q}}
\newcommand{\Z}{\mathbf{Z}}

\newcommand{\y}{\mathbf{y}}
\newcommand{\f}{\mathbf{f}}

\newcommand{\U}{\mathbf{u}}

\begin{document}

\title{A self adjusting multirate algorithm based on the TR-BDF2  method}
\author{Luca Bonaventura$^{(1)}$, Francesco Casella$^{(2)}$ \\
Ludovica Delpopolo$^{(1)}$, Akshay Ranade$^{(3)}$ }
\maketitle

\begin{center}
{\small
$^{(1)}$ MOX -- Modelling and Scientific Computing \\
Dipartimento di Matematica, Politecnico di Milano \\
Via Bonardi 9, 20133 Milano, Italy\\
 {\tt luca.bonaventura@polimi.it}, {\tt ludovica.delpopolo@polimi.it}
}
\end{center}

\begin{center}
{\small
$^{(2)}$  Dipartimento Elettronica, Informazione e Bioingegneria\\
 Politecnico di Milano \\
Via Ponzio 34/35, 20133 Milano, Italy\\
{\tt francesco.casella@polimi.it}
}
\end{center}

\begin{center}
{\small
$^{(3)}$  Department of Computer Science\\
University College Cork\\
Western Road\\
Cork, Ireland\\
{\tt akshay.ranade@cs.ucc.ie}
}
\end{center}

\date{}

\noindent
{\bf Keywords}:  Multirate methods, stiff ODE problems, Implicit Runge Kutta methods, SDIRK methods

\vspace*{1.0cm}

\noindent
{\bf AMS Subject Classification}:  65L04, 65L05, 65L07, 65M12, 65M20

\vspace*{1.0cm}

\pagebreak

\abstract{We propose a self adjusting multirate method based on the TR-BDF2 solver. The potential advantages
of using  TR-BDF2  as the key component of a 
 multirate framework are highlighted. A linear stability analysis of the resulting approach is 
presented and the stability features of the resulting  algorithm are analysed. 
The analysis framework is completely general and allows to study along the same lines the stability of self adjusting multirate methods based on a generic one step solver. A number of numerical experiments demonstrate the efficiency and accuracy of the resulting approach also the time discretization of hyperbolic partial differential equations.}

\pagebreak

\section{Introduction}
 \label{intro} \indent

 The concept of  multirate methods was first proposed in \cite{rice:1960} and a great attention
 has been devoted to these methods, see e.g. \cite{andrus:1979}, \cite{gear:1984} and the more comprehensive review  reported in \cite{ranade:2016}.  The basic idea of multirate methods is to integrate each component using a time step that is the most appropriate to the timescale of that component. Slowly varying components are integrated with larger time steps, while smaller time steps are used for fast components only. Thus, multirate methods can avoid a significant amount of the computation that is necessary in the standard, single rate approach. In other words, in the multirate approach,  the most appropriate time resolution is employed for each  state variable of the system. In the context of multirate algorithms, the faster components  are often referred to as  active and the slow ones are called instead  latent. 
  It is also to be remarked that conceptually similar approaches
 have been introduced   in the  literature on numerical weather prediction, see e.g. \cite{klemp:1978},
  for the purpose of avoiding the time step stability restrictions related  to fast propagating acoustic waves.

 The specific way in which the partitioning is performed and the degree of coupling
 allowed between latent and active components  determines the efficiency and the stability of the resulting methods. The stability of  multirate methods has been investigated, among others,  in \cite{andrus:1993} and \cite{gunther:2001}, while the stability  of split-explicit methods has been analyzed in   \cite{baldauf:2010}, \cite{skamarock:1992}.  
  In general, it is seen that the stability region of the multirate methods decreases as the strength of the coupling between the partitions  or the stiffness of the system increases. 
   In spite of this, there is ample evidence in the literature that multirate algorithms can offer a
    promising path for the numerical integration of large systems with widely varying time scales 
    or with highly localised periods of activity.
 In order to address the potential problems of multirate methods,
 the so-called compound step multirate methods have been proposed in  \cite{gunther:2001}. 
A self-adjusting, recursive time stepping strategy has been   proposed in  \cite{savcenco:2007}.
In this kind of approach, a tentative  global  step is first taken for all components, using a robust, unconditionally stable method. The time step is then reduced
only  for those components for which some local error estimator is greater than the specified tolerance.
In this way, automatic detection of the latent components is achieved.
 In \cite{savcenco:2008}, an algorithm based on this recursive strategy was shown to be more stable than the original compound step algorithm, while a stability analysis of a self-adjusting  algorithm based on the $\theta-$method  has been proposed in  \cite{hundsdorfer:2009}.
 
 In this paper, we propose a self-adjusting multirate method based on the same concept as
 \cite{savcenco:2007}, that relies on the TR-BDF2 method as fundamental single rate solver.  The TR-BDF2 method has been originally introduced in \cite{bank:1985} and more thoroughly analyzed in  \cite{hosea:1996}. 
 It is a second order, one step,
 L-stable implicit method endowed with a number  of interesting properties, as discussed in  \cite{hosea:1996}.
  For example, it allows for cheap error estimation
 via an embedded third order method, it can provide  continuous output via a cubic Hermite interpolant  and 
 it has an explicit second order
 companion with which it can be combined to form a second order, implicit-explicit additive Runge Kutta method.
 Due to its favourable properties, it has been recently applied for efficient discretization
 of high order finite element methods for numerical weather forecasting in \cite{giraldo:2013}, \cite{tumolo:2015}.
 Its appealing features  also make it   an interesting candidate for multirate approaches targeted at  large scale stiff problems, resulting either from  the simulation industrial systems or from the spatial discretization of partial differential equations.
 In particular, the  cubic Hermite interpolant can be employed to maintain higher accuracy also for stiff
 problems in the framework of a multirate approach.
  The proposed self-adjusting multirate TR-BDF2 method   is  equipped here with a partitioning 
  and time step selection criterion based on the technique proposed in  \cite{fok:2015}.  
  
  The stability of the method has been analyzed in the framework of the classical linear model problem.
  Even though no complete theoretical results can be achieved, numerical computation of the stability function norm
  for a range of relevant model problems shows that the method has  good stability properties.
  This is true   not only for contractive problems with strictly dissipative behaviour, but also for problems
  with purely  imaginary eigenvalues, which suggests that the method could be advantageous also for
  the application to hyperbolic PDEs and structural mechanics problems.
  The numerical results obtained on several benchmarks  show that the application of the proposed method leads
  to significant efficiency gains, that are analogous to those achieved by other self-adjusting multirate approaches
  with respect to their single rate counterparts. In particular, the method appears to perform well when applied to the time discretization of nonlinear, hyperbolic partial differential equations, allowing to achieve automatic detection
  of complex localized phenomena such as shock waves and significant reductions in computational cost.

 In section \ref{selfa}, the self-adjusting implicit multirate approach is described.
    The time step refinement criterion and the resulting
   partitioning principle between active and latent components are presented  in section \ref{refinement}. 
 The TR-BDF2 method is reviewed in section \ref{solvers}.
 The interpolation procedures we have considered are described in section
 \ref{interpolation}. A linear stability analysis is then presented in section \ref{stability}.
 The approach proposed for such analysis is rather general and allows to study generic multirate methods based on
 any one step solver and interpolation procedure. 
 Numerical simulations are presented in section \ref{tests}, while some conclusions
 and perspectives for future work are presented in \ref{conclu}.

  \section{A self adjusting implicit multirate approach}
 \label{selfa} 
 \indent
 
 We focus on multirate methods for the solution of the Cauchy initial value problem
\begin{equation}\label{eq:cauchy}
\mathbf{y}^{\prime}=\mathbf{f}(t,\mathbf{y}) \ \ \ \ \ \ \mathbf{y}(0)=\mathbf{y}_0\in\mathbf{R}^m, \ \ \ \ \ t\in[0,T].
\end{equation}
We   consider time discretizations associated to discrete time levels
$t_n,\  n=0,\dots,N $ such that $h_n=t_{n+1}-t_n$ and
we will denote by
$\mathbf{u}_n $ the numerical approximation of $\mathbf{y}(t^n).$ 
 We will also denote by ${\bf u}^{n+1}={\cal S} ({\bf u}^{n},h_n) $ the
implicitly defined operator ${\cal S} :\mathbf{R}^m\rightarrow \mathbf{R}^m $
 whose application is equivalent to the computation of  one step of size $h_n$  of
a given single step method. While here only implicit methods will be considered,
the whole framework can be extended to explicit and IMEX methods.
Notice that, if $\mathbf{P}$ is the projector onto a linear subspace ${\cal V} \subset {\mathbf{R}}^m $ with dimension $p < m,$
the operator ${\cal S}^{\cal V} :\mathbf{R}^p\times \mathbf{R}^{m-p}\rightarrow \mathbf{R}^p $
that represents the solution of the subsystem obtained freezing the components
of the unknown belonging to  ${\cal V}^{\perp}  $  to the value ${\bf z} \in \mathbf{R}^{m-p}$ can be defined by
${\bf y} =   {\cal S}^{\cal V} ({\bf x},{\bf z}, h_n )= \mathbf{P}{\cal S} ( {\bf x}\oplus\mathbf{z}, h_n) .$
Furthermore, we will denote by $\mathbf{Q} ({\bf u}^{n+1},{\bf u}^{n},\zeta) $ 
an interpolation operator, that provides an approximation of the numerical solution at
intermediate time levels $t_n+\zeta,$ where $\zeta \in[0,h_n].$ Linear interpolation is often employed, but, for
 multistage methods, knowledge of the intermediate stages  also allows the application of more accurate interpolation procedures  without 
substantially increasing the computational cost.

We propose a self-adjusting  multirate algorithm that
  is a generalization of the method introduced in \cite{savcenco:2007}, 
in which the choice of the local time step is left completely to the time step selection criterion.  
The multirate algorithm can be described as follows.
\begin{itemize}
	\item[1)]  Perform  a tentative global (or macro) time step of size $h_n $ with the standard single rate method and compute
	 $\hat{\bf u}^{n+1}= {\cal S} ({\bf u}^{n},h_n).$ The specific way in which error control mechanisms are
	  applied to the choice of the global and local time steps
	are crucial for the efficiency of the algorithm and will be discussed in detail in section
	\ref{refinement}.  
	 \item[2)]  Compute the error estimator 
	 and  partition the state space into  active and  latent variables, based on the value of the error estimate.
	  The projection onto the subspace ${\cal V}_{0}$
	of the active variables is denoted by $\mathbf{P}_n=\mathbf{P}_n^{(0)},$  
	while the projection onto the complementary
	subspace will be denoted by $\bar{\mathbf{P}}_n^{(0)}.$
	Define $    \bar{\mathbf{P}}_n^{(0)}{\bf u}^{n+1} =    \bar{\mathbf{P}}_n^{(0)}\hat {\bf u}^{n+1} $
	as well as $ {\bf u}^{n,0}= {\bf u}^{n} $ and
	  $t_{n,0}=t_{n}. $ 
	\item[3)]  For $k\geq 1,$
	choose a local (or micro) time step  $h_n^{(k-1)} $ for the active variables, based on the value of the error estimator.  Set  $t_{n,k}=\min\{t_{n,k-1}+h_n^{(k-1)},  t_{n+1}\}. $ 
	\begin{itemize}
		\item[3.1)]      Update the latent variables by interpolation
	$$	
	 \bar{\mathbf{P}}_n^{(k-1)}{\bf u}^{n,k}=
	\mathbf{Q}  ( \bar{\mathbf{P}}_n^{(k-1)}{\bf u}^{n+1},  \bar{\mathbf{P}}_n^{(k-1)} {\bf u}^{n,k-1}, h_n^{(k-1)} ).
	$$
	\item[3.2)]  Update the active variables  by computing
		$$
	 \mathbf{P}^{(k-1)}_n{\bf u}^{n,k}=
	  {\cal S}^{{\cal V}_{k-1}} (\mathbf{P}^{(k-1)}_n{\bf u}^{n,k-1},  \bar{\mathbf{P}}_n^{(k-1)}{\bf u}^{n,k}, h_n^{(k-1)}).
	$$
	\item[3.3)] Compute the error estimator  for the active variables only and   partition again  $ {\cal V}_{k-1} $
	into latent and active variables. Denote by ${\cal V}_{k}\subset {\cal V}_{k-1} $ the new
	subspace of active variables and  by  $\mathbf{P}_n^{(k)} $  the corresponding projection.
	
	\item[3.4)] Repeat  3.1) - 3.3) until $t_{n,k}=t_{n+1}.$
%
		\end{itemize}
	\end{itemize}

Clearly, the effectiveness of the above procedure depends in a crucial way on the accuracy and stability of the basic ODE 
solver ${\cal S}, $ as well as on the time step
refinement and partitioning criterion.  Furthermore, as well known in multirate methods,
unconditional stability of the solver ${\cal S} $ does not necessarily entail that the same property holds for
the derived multirate solver.
The refinement and partitioning criterion will be described in detail in section
\ref{refinement}. Here, we only remark that, as it will be shown by the numerical experiments
in section \ref{tests}, the approach outlined above is able to reduce significantly the computational
cost with respect to the equivalent single rate methods, without a major reduction in stability
for most of the envisaged applications.

   \section{The time step refinement and partitioning criterion}
 \label{refinement} 
 \indent
 
 We now describe in detail the time step refinement and partitioning strategy that we have used  in the multirate algorithm described in section \ref{selfa}.
Our approach is based on the strategy proposed for an explicit Runge Kutta multirate method in
 \cite{fok:2015}, where the time steps for refinement are obtained from the error estimates of the global step. The user specified tolerance plays an important role in the partitioning of the system. In \cite{fok:2015}, a simple absolute error tolerance was considered. However, in most engineering system simulations, whenever the typical values of different components can vary greatly, the tolerance used is in general
 a combination  of absolute and relative error tolerances, see e.g. \cite{soderlind:2006}.
  We have thus extended the strategy  proposed in
 \cite{fok:2015} in order to employ such a combination of absolute and relative error tolerances.
 
More specifically, denote by   $\tau_r, $ $\tau_a$ the user defined error tolerances for relative and absolute errors, respectively.
Furthermore, assume that the tentative global step $\hat{\bf u}^{n+1}={\cal S} ({\bf u}^{n},h_n^{(0)}) $
has been computed and that   an   error estimator $\bm \epsilon^{n+1,0}  $ is available.
The first task of the time step selection criterion is to asses whether the global time step $h_n^{(0)} $ has been properly chosen.
Denoting by $\epsilon^{n+1,0}_i, \hat u^{n+1}_i, i=1, \dots,m $ the $i-$th components of the error estimator and of the
tentative global step numerical solution, respectively, 
we define a   vector $\bm \eta^{n+1,0} $ of normalized errors with components
   $$ \eta^{n+1,0}_i= \frac{\epsilon^{n+1,0}_i}{\tau_r|   \hat u^{n+1}_i| + \tau_a}. $$
   Since clearly the condition $ \max_i  \eta^{n+1,0}_i \leq 1$
has to be enforced, before proceeding to the partitioning into
active and latent variables, this condition is checked and the global step
is repeated with a smaller value of $h_n^{(0)} $ whenever it is violated.
Numerical experiments have shown that, while an increase of efficiency with respect
to the single rate version of the method is always achieved, independently of the choice of the tentative
step,  the greatest improvements are only possible if the global time step  does not have
to  be repeated too often.

   Once the condition  $ \max_i  \eta^{n+1,0}_i \leq 1$ is satisfied,
   the set of indices of the components flagged for the first level of time step refinement is identified by 

\begin{equation}
	{\cal A}^{n+1,0} = \{ i :  \eta^{n+1,0}_i > \delta \| \bm \eta^{n+1,0} \|_{\infty}\},
\end{equation}
where $\delta\in(0,1)$ is a user defined coefficient. The smaller the value of $\delta, $
the larger is the fraction of  components marked for refinement. Notice that,
if $\delta$ is set to unity,  the algorithm effectively operates in  single rate mode.
For each iteration $k=1,\dots,k_{max} $ of the algorithm described in section \ref{selfa},
active variable sets ${\cal A}^{n+1,k} $ are then defined analogously.
In each iteration $k=0,\dots,k_{max}, $ the time step $h_n^{(k)} $ is chosen according
to the standard criterion (see e.g. \cite{lambert:1991}, \cite{savcenco:2007}) 

\begin{equation}
\label{eq:refinement}
	h_n^{(k)} = \nu \min_{j\in{\cal A}^{n+1,k} }
	\left( \dfrac{\tau_r|\hat u^{n+1}|_j +\tau_a}{\epsilon^{n+1,k}_i }\right)^{\frac{1}{p+1}}, \hspace{2em}
\end{equation}
where $p $ is the convergence order of the solver ${\cal S} $ and $\nu $ is an user defined 
safety parameter.
 
 \section{The single rate TR-BDF2 method}
 \label{solvers} 
 \indent
 
 The TR-BDF2 method was originally proposed in \cite{bank:1985}. It is  a composite one step, two stage method, consisting
of one stage of the trapezoidal method followed by another of the BDF2 method. The stages are so adjusted that both the trapezoidal and the BDF-2 stages use the same Jacobian matrix. This composite method  has been reinterpreted 
in \cite{hosea:1996} 
as a one step Diagonally Implicit Runge Kutta (DIRK) method with two internal stages  The TR-BDF2 method is also in some sense the optimal method among all 3-stage DIRK methods,  owing to the following properties:
\begin{itemize}
	\item[1)] it is first same as last (FSAL), so that only two implicit stages need to be evaluated; 
	\item[2)] the Jacobian matrix for both the stages is the same;
	\item[3)] it has an embedded third order companion that allows for a cheap error estimator;
	\item[4)] the method is strongly S-Stable;
	\item[5)] all the stages are evaluated within the time step;
	\item[6)] it is endowed with a cubic Hermite interpolation algorithm for dense output that yields globally $C^1$ continuous
	trajectories.
\end{itemize}
Features 3), 5) and 6) will play an important role in the multirate method proposed here.
Furthermore, as shown in  \cite{giraldo:2013}, the method   has an explicit second order
 companion with which it can be combined to form a second order implicit-explicit additive Runge Kutta method.
 Due to its favourable properties, it has been recently applied for efficient discretization
 of high order finite element methods for numerical weather forecasting in \cite{giraldo:2013}, \cite{tumolo:2015}.

The TR-BDF2 method, considered as a composite method consisting of a step with the trapezoidal method followed by a step of the BDF2 method, can be written as

\begin{equation}
\begin{aligned}
	{\bf u}^{n+\gamma}&={\bf u}^n + \frac{\gamma h_n}2\left ({\bf f}(t,{\bf u}^n) + {\bf f}(t,{\bf u}^{n+\gamma}) \right)\\
	{\bf u}^{n+1}&= \frac{1}{\gamma(2-\gamma)}{\bf u}^{n+\gamma} -\frac{(1-\gamma)^2}{\gamma(2-\gamma)}{\bf u}^n + \frac{1-\gamma}{2-\gamma}h_n
	{\bf f}(t,{\bf u}^{n+\gamma})
\end{aligned}
\label{trbdf2_comp}
\end{equation}

The TR-BDF2 method viewed as a DIRK method has the following Butcher tableau:

\vspace{4em}
\setlength{\tabcolsep}{1.5em}
\begin{center}
\begin{math}
\begin{tabular}{r|rrr}
$0$ & ~$0$ & ~$0$ & ~$0$\\
$\gamma$& ~$d$ & ~$d$ & ~$0$\\
$1$& ~$w$ & ~$w$ & ~$d$\\\hline
~& ~$w$ & ~$w$ & ~$d$\\\hline
~*& $\frac{(1-w)}{3}$ & $\frac{(3w+1)}{3}$& $\frac{~d~}{3}$
	\end{tabular}
	\label{tab:butchertrbdf2}
\end{math}
\end{center}
where $\gamma=2-\sqrt{2}$, $d=\frac{\gamma}{2}$ and $w=\frac{\sqrt{2}}{4}$ and the row * corresponds to the embedded
 third order method that can be used to build a convenient  error estimator.
 Here, the value  $\gamma=2-\sqrt{2} $ is chosen for the two stages to have the same Jacobian matrix, 
 as well as in order to achieve L-stability, as it will be shown shortly.
It can be seen that the method has the First Same As Last (FSAL) property, i.e., the first stage of any step is the same as the last stage of the previous step. Thus, in any step, the first explicit stage need not be computed. 

 The implementation of the two implicit stages is done as suggested in \cite{hosea:1996}.
  The two stages according to the Butcher tableau are given by equation
\eqref{trbdf2_comp}.
Instead of iterating on the variable ${\bf u}^{n+1}$, 
we define another variable ${\bf z}  =h{\bf f}(t,{\bf u})$ and solve for this variable by iteration. Thus, for the first implicit stage, we take \mbox{${\bf u}^{n+\gamma,k}={\bf u}^n + d{\bf z}^n +d{\bf z}^{n+\gamma,k}$} and ${\bf z}^{n+\gamma,k}$ is computed by Newton iterations as

\begin{equation*}
\begin{aligned}
	(\mathbf{I}-dh\mathbf{J}^n)\bm \Delta^k&=h_n{\bf f}(t_{n+\gamma},\y_{n+\gamma}^k) - {\bf z}^{n+\gamma,k}\\
	{\bf z}^{n+\gamma,k+1} &= {\bf z}^{n+\gamma} + \bm \Delta^k,
\end{aligned}
\end{equation*}
where  ${\bf J}^n= {\bf J}(t_n, {\bf u}_n)$ denotes the Jacobian of ${\bf f}.$  
Similarly, for the second implicit stage we take ${\bf u}_{n+1}^k = {\bf u}_n + w{\bf z}^n + w {\bf z}^{n+\gamma} + d{\bf z}^{n+1,k}$ and the Newton iteration is given by
\begin{equation*}
\begin{aligned}
	(\mathbf{I}-dh\mathbf{J}^n)\bm \Delta^k&=h{\bf f}(t_{n+\gamma},\y_{n+\gamma}^k) - {\bf z}^{n+\gamma,k}\\
	{\bf z}^{n+1,k+1} &= {\bf z}^{n+1} + \bm \Delta^k.
\end{aligned}
\end{equation*}

The reason for doing so is that if the problem is stiff, a function evaluation will amplify the numerical error in the stiff components. Iterating on ${\bf z}$, however, ensures that it is computed to within the specified tolerance. For a more detailed discussion the interested reader is referred to \cite{hosea:1996}.
An important point to be noted is that, although the TR-BDF2 method is L-stable, its third order companion formula is not. Therefore, the error estimate at time level $n+1,$   given by
\begin{equation*}
	\bm\epsilon^{n+1,*}=(b_1^* - b_1){\bf z}^n +(b_2^* - b_2){\bf z}_{n+\gamma} +(b_3^* - b_3){\bf z}_{n+1} 
\end{equation*}
cannot be expected to be accurate for stiff problems. To overcome this problem, it is suggested
in  \cite{hosea:1996}   
to modify the error estimate by considering as error estimator the quantity $\bm \epsilon_{n+1} $ defined as the solution
of the linear system
 
\begin{equation}
\label{mod_errest}
	(\mathbf{I}-dh\mathbf{J}^n)\bm \epsilon^{n+1}=\bm \epsilon^{n+1,*}.
\end{equation}
This modification of the error estimate  allows to improve it for stiff components, while preserving its accuracy in the limit of small time steps. Notice that the stability function of the TR-BDF2 method is given by
\begin{equation}
 \label{eq:trbdf2_stability_function}
 R(z) =\frac{[1 + (1 - \gamma)^2] z + 2(2-\gamma)}
      {z^2(1-\gamma)\gamma + z(\gamma^2-2) +2(2-\gamma)},
\end{equation}
whence it can be seen that the method is L-stable for $\gamma = 2-\sqrt{2}.$
 
%

   \section{The interpolation procedures}
 \label{interpolation} 
 \indent
 
 An essential component of any multirate algorithm is the  procedure employed
 to reconstruct the values of the latent variables at those intermediate time levels
 for which only the active variables are computed. Self-adjusting multirate procedures
 based on implicit methods, such as the one proposed in \cite{savcenco:2007} and that
 presented in this paper, can avoid the use of extrapolation, thus increasing their overall stability.
 The simplest and most commonly used procedure is linear interpolation, that is defined for $\zeta\in[0,h_n] $ as
 \begin{equation}
 \label{lin_int}
 \mathbf{Q} ({\bf u}^{n+1},{\bf u}^{n},\zeta) = \frac{\zeta}{h_n} {\bf u}^{n+1} + \frac{(h_n-\zeta)}{h_n} {\bf u}^{n}.
 \end{equation}
 In case a multi-stage solver is employed, higher order interpolation can also be employed at no extra cost.
 Assuming that a sufficiently accurate approximation of the solution  ${\bf u}^{n+\lambda} $ is available at time $t_{n+\lambda}\in (t_n,t_{n+1}) $
 and setting $h_{\lambda}= t_{n+\lambda}-t_n,$
 quadratic Lagrange interpolation can then be written as 
 \begin{eqnarray}
 \label{quad_int}
 \mathbf{Q} ({\bf u}^{n+1},{\bf u}^{n},\zeta) &=&\frac{(\zeta-h_{\lambda})\zeta }{h_{n}(h_{n}-h_{\lambda})}{\bf u}^{n+1}  
+ \frac{\zeta (\zeta -h_n) }{(h_{\lambda}-h_n)h_{\lambda}}{\bf u}^{n+\lambda}\\
 &+&\frac{(h_{\lambda}-\zeta)(h_{n}-\zeta) }{h_{\lambda}h_{n}} {\bf u}^{n}. \nonumber
 \end{eqnarray}
 
 An interesting feature of the TR-BDF2 method in the multirate framework is that the method, as explained in \cite{hosea:1996},
is endowed with a cubic Hermite, globally $C^1$ interpolation algorithm for dense output. 
Using the notation of section   \ref{solvers}, the Hermite cubic interpolant can be defined as
\begin{eqnarray}
 \mathbf{Q} ({\bf u}^{n+1},{\bf u}^{n},\zeta) &=& (\bm \alpha_3 - 2\bm \alpha_2)\beta^3(\zeta) \\
 &+& (3\bm\alpha_2 - \bm\alpha_3)\beta^2(\zeta)   
  + \bm\alpha_1\beta(\zeta) +\bm\alpha_0,  \nonumber
 \label{hermite}
\end{eqnarray}
	where, for $ 0\leq \zeta \leq \gamma h_n $ one has
	\begin{equation}
\begin{aligned}
	\bm\alpha_0&=\U_n, \hspace{1em} \bm\alpha_1=\gamma \mathbf{z}_n, \hspace{1em} \bm\alpha_2=\U_{n+\gamma}-\U_n -\gamma\bm z_n, \nonumber \\
	\bm\alpha_3&=\gamma(\bm z_{n+\gamma} - \bm z_n), \hspace{1em}   \hspace{1em} \beta(\zeta)= \frac{\zeta}{\gamma h_n} \nonumber
	 \end{aligned}
\end{equation}
and  for $\gamma h_n\leq \zeta\leq h_{n}  $ one has instead
\begin{equation}
\begin{aligned}	
	\bm\alpha_0&=\U_{n+\gamma}, \hspace{0.25em}\bm\alpha_1=(1-\gamma) \mathbf{z}_{n+\gamma}, \hspace{0.25em} \bm\alpha_2=\U_{n+1} - \U_{n+\gamma} - \bm\alpha_1,  \nonumber\\
	\bm\alpha_3&=(1-\gamma)( \mathbf{z}_{n+1} - \mathbf{z}_{n+\gamma}) \hspace{0.75em}  \hspace{0.75em} \beta(\zeta)=\frac{\zeta-\gamma h_n}{(1-\gamma)h_n}.  \nonumber
\end{aligned}
\end{equation}
This interpolant allows to have accurate
approximations of the latent variables without extra computational cost or memory storage requirements, since  
it employs the intermediate stages of the TR-BDF2 method in order to achieve higher order accuracy. This contributes
to making  the TR-BDF2 method an attractive choice as the basis for a multirate approach.

  \section{Linear stability analysis}
 \label{stability} 
 \indent
 
 In this section,  we will study the linear stability of the multirate algorithm outlined in the
 previous sections in the case of a generic  one step  solver.
Following the approach employed in \cite{hundsdorfer:2009} for the $\theta-$method,
 we will consider the special case in which the
macro time step is constant and denoted by $h $ and a single refinement is carried out, with the
refined time step used for the active components given by $h/2.$
As a result, a single subspace of active variables is introduced,
${\cal V} \subset {\mathbf{R}}^m, $ with dimension $p < m.$ 
As in section \ref{selfa}, $\mathbf{P}$ will denote the projector onto  this subspace, while
$\mathbf{P}^{\perp}$ will denote the corresponding projector onto the complementary
subspace ${\cal V}^{\perp}  $  of the latent variables. We will also denote by 
$\mathbf{E}, \mathbf{E}^{\perp}$ the corresponding natural embedding operators
of ${\cal V},  {\cal V}^{\perp}  $ into $ {\mathbf{R}}^m. $ As a result,
the operators  $\boldsymbol{\Pi}=\mathbf{E}\mathbf{P} $ and $\boldsymbol{\Pi}^{\perp}=\mathbf{E}^{\perp}\mathbf{P}^{\perp} $ 
represent the zeroing of the components of a given vector belonging to  ${\cal V},  {\cal V}^{\perp}  $
respectively.  It is to be remarked that this special setting is a major simplification with respect to the generality of the algorithm
described in section \ref{selfa}, so that the results of the analysis may not fully explain the behaviour of the algorithm
in practical applications. In particular, large values of the stability function norm in this context do not necessarily imply that
the method is unable to achieve the desired accuracy when applied in its most general form.

We  now consider the linear problem 	$\y'=\f(t,\y)=\mathbf{A}\y$  and we set $h\mathbf{A}=\Z$. 
We assume that the amplification matrix of  the basic one step solver can be written  as
$$ \mathbf{R}= \mathbf{R}( \mathbf{Z})= \mathbf{D}( \mathbf{Z})^{-1} \mathbf{N}( \mathbf{Z})=\mathbf{D}^{-1} \mathbf{N}, $$
where $\mathbf{D}, \mathbf{N}$ are polynomials in $\mathbf{Z}.$ Notice that this includes most standard one step
solvers, either explicit or implicit. We will also use the notations
$$   \mathbf{D}_{1/2}  =  \mathbf{D}( \mathbf{Z}/2 ) \ \ \ \  \mathbf{N}_{1/2}=\mathbf{N} ( \mathbf{Z}/2). 
 $$

Notice that, in the special case considered here, the interpolation operator can be represented by the application to the data
of a matrix $\Q_{1/2} $ whose entries only depend on $h $  and on the amplification matrix of the single rate method employed.
For example, the linear interpolant \eqref{lin_int} is represented in this case by the matrix
$$\Q_{1/2} = \frac 12 \left (  \mathbf{I}+\mathbf{R} \right ), $$ 
while for the cubic Hermite interpolant   \eqref{hermite} one has
$$\mathbf{Q}_{1/2} =   \mathbf{I} +  \beta \gamma \mathbf{Z}+   \beta^2\mathbf{F}  
+ \beta^3  \mathbf{G}, $$
where $\beta=1/(2\gamma) $ and
$$  \mathbf{F}=3 (\mathbf{R}_{\gamma}  - \mathbf{I} -{\gamma}\mathbf{Z})  -\gamma \mathbf{Z} (\mathbf{R}_{\gamma}  - \mathbf{I})
 \ \ \ \ \   \mathbf{G}=\gamma \mathbf{Z} (\mathbf{R}_{\gamma}  - \mathbf{I})
 -2(\mathbf{R}_{\gamma}  - \mathbf{I} -{\gamma}\mathbf{Z}),     $$
$$ \mathbf{R}_{\gamma}= \left (  \mathbf{I} -\frac{\gamma}2 \mathbf{Z}\right )^{-1} 
\left ( \mathbf{I} +\frac{\gamma}2 \mathbf{Z} \right ).$$
 
 Given these definitions, one step of 
 the self-adjusting multirate algorithm can be rewritten in this special setting as

\begin{eqnarray}
		\label{linm1}	
	&&	\mathbf{P}^{\perp}{{\bf u}}^{n+1}= \mathbf{P}^{\perp}\mathbf{R}{\bf u}^{n}=\mathbf{P}^{\perp}{\bf u}^{n,2}
	 \ \ \ \ \	\mathbf{P}^{\perp}{\bf u}^{n,1}=\mathbf{P}^{\perp}\Q_{1/2}{\bf u}^n 	\nonumber \\
	&& \mathbf{P} \mathbf{D}_{1/2} \left (\boldsymbol{\Pi}{\bf u}^{n,1}+ \boldsymbol{\Pi}^{\perp}\Q_{1/2}{\bf u}^n
	\right )=  \mathbf{P}\mathbf{N}_{1/2}{\bf u}^n 	\\
	&&\mathbf{P} \mathbf{D}_{1/2} \left (\boldsymbol{\Pi}{\bf u}^{n,2}+\boldsymbol{\Pi}^{\perp} \mathbf{R}{\bf u}^{n}  \right ) =  \mathbf{P}\mathbf{N}_{1/2}(\boldsymbol{\Pi}{\bf u}^{n,1}+\boldsymbol{\Pi}^{\perp}\Q_{1/2}{\bf u}^n ).	\nonumber
		\end{eqnarray}
Recalling the previous definitions, the equation for 	$\mathbf{P}{\bf u}^{n,1}	$ can be rewritten as
$$
\mathbf{P} \mathbf{D}_{1/2}\mathbf{E}\mathbf{P}{\bf u}^{n,1}= \mathbf{P} \mathbf{N}_{1/2}{\bf u}^{n}
-\mathbf{P} \mathbf{D}_{1/2}\boldsymbol{\Pi}^{\perp}\Q_{1/2}{\bf u}^n,
$$
which implies	
$$
\mathbf{P}{\bf u}^{n,1}= (\mathbf{P} \mathbf{D}_{1/2}\mathbf{E})^{-1}\left [\mathbf{P} \mathbf{N}_{1/2} 
-\mathbf{P} \mathbf{D}_{1/2}\boldsymbol{\Pi}^{\perp}\Q_{1/2} \right] {\bf u}^n.
$$
Substituting this expression in the last equation of \eqref{linm1} one obtains
\begin{eqnarray}
		\label{linm2}
	&&\mathbf{P} \mathbf{D}_{1/2} \mathbf{E}\mathbf{P}{\bf u}^{n,2}
	 =    - \mathbf{P} \mathbf{D}_{1/2}  \boldsymbol{\Pi}^{\perp}\mathbf{R}{\bf u}^{n} \nonumber \\
	&&  +\mathbf{P}\mathbf{N}_{1/2}\mathbf{E}   
	(\mathbf{P} \mathbf{D}_{1/2}\mathbf{E})^{-1}
\left [\mathbf{P} \mathbf{N}_{1/2} -\mathbf{P} \mathbf{D}_{1/2}\boldsymbol{\Pi}^{\perp}\Q_{1/2} \right] {\bf u}^n \\
	  &&+\mathbf{P}\mathbf{N}_{1/2}\boldsymbol{\Pi}^{\perp}\Q_{1/2}{\bf u}^n .	\nonumber
		\end{eqnarray}	
The stability matrix of the complete multirate method	can then be obtained	
deriving from \eqref{linm2} an explicit expression for
$\mathbf{P}{\bf u}^{n,2}$ and substituting in
$${\bf u}^{n+1}=	\boldsymbol{\Pi}{\bf u}^{n,2} +\boldsymbol{\Pi}^{\perp}\mathbf{R}{\bf u}^{n}. $$
Introducing  for compactness the notation
$$
 \mathbf{A}^{\parallel,\parallel}= \mathbf{P} \mathbf{A}\mathbf{E}, \ \ \mathbf{A}^{\parallel,\perp}= \mathbf{P} \mathbf{A}\mathbf{E}^{\perp}
$$
one can express the stability matrix as follows 
  \begin{eqnarray}
		\label{stabmat}
\mathbf{R}_{mr}&=&\mathbf{E}\left \{ \left (\mathbf{D}_{1/2}^{\parallel,\parallel}\right)^{-1}\mathbf{N}_{1/2}^{\parallel,\parallel}
 \left (\mathbf{D}_{1/2}^{\parallel,\parallel}\right)^{-1} \left [   \mathbf{P} \mathbf{N}_{1/2}
-  \mathbf{D}_{1/2}^{\parallel,\perp}\mathbf{P}^{\perp}\Q_{1/2}  \right] \right. \nonumber \\
&+& \left. \left (\mathbf{D}_{1/2}^{\parallel,\parallel}\right)^{-1} \mathbf{N}_{1/2}^{\parallel,\perp}\mathbf{P}^{\perp}\Q_{1/2} 
- \left (\mathbf{D}_{1/2}^{\parallel,\parallel}\right)^{-1}\mathbf{D}_{1/2}^{\parallel,\perp}\mathbf{P}^{\perp} \mathbf{R}
\right\} +\boldsymbol{\Pi}^{\perp}\mathbf{R}. \nonumber
\end{eqnarray}		
 
 We will not  attempt a complete theoretical analysis based on the previously derived expression of
 $\mathbf{R}_{mr}.$ However, the same expression can be employed to compute numerically the stability matrix norm and
 spectral radius for relevant model problems. We have done so for three classes of linear ODE systems $\y'=\mathbf{A}\y.$ 
 Notice that, for simplicity, all the systems have been rearranged so that the latent variables are the first
 of the vector of unknowns.
 
\begin{figure}[!htbp]
\centering
\includegraphics[width=0.46\textwidth]{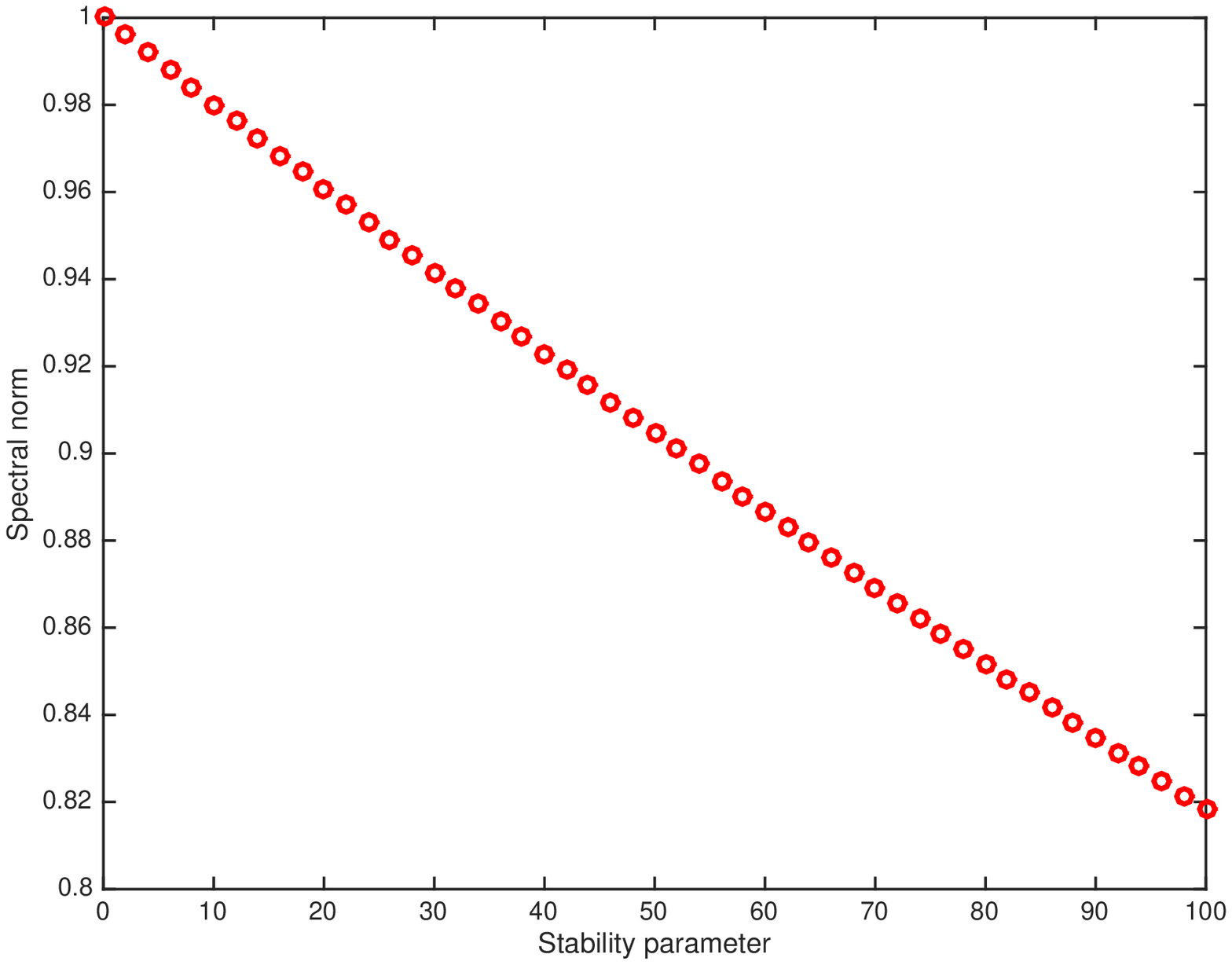}a)
\includegraphics[width=0.46\textwidth]{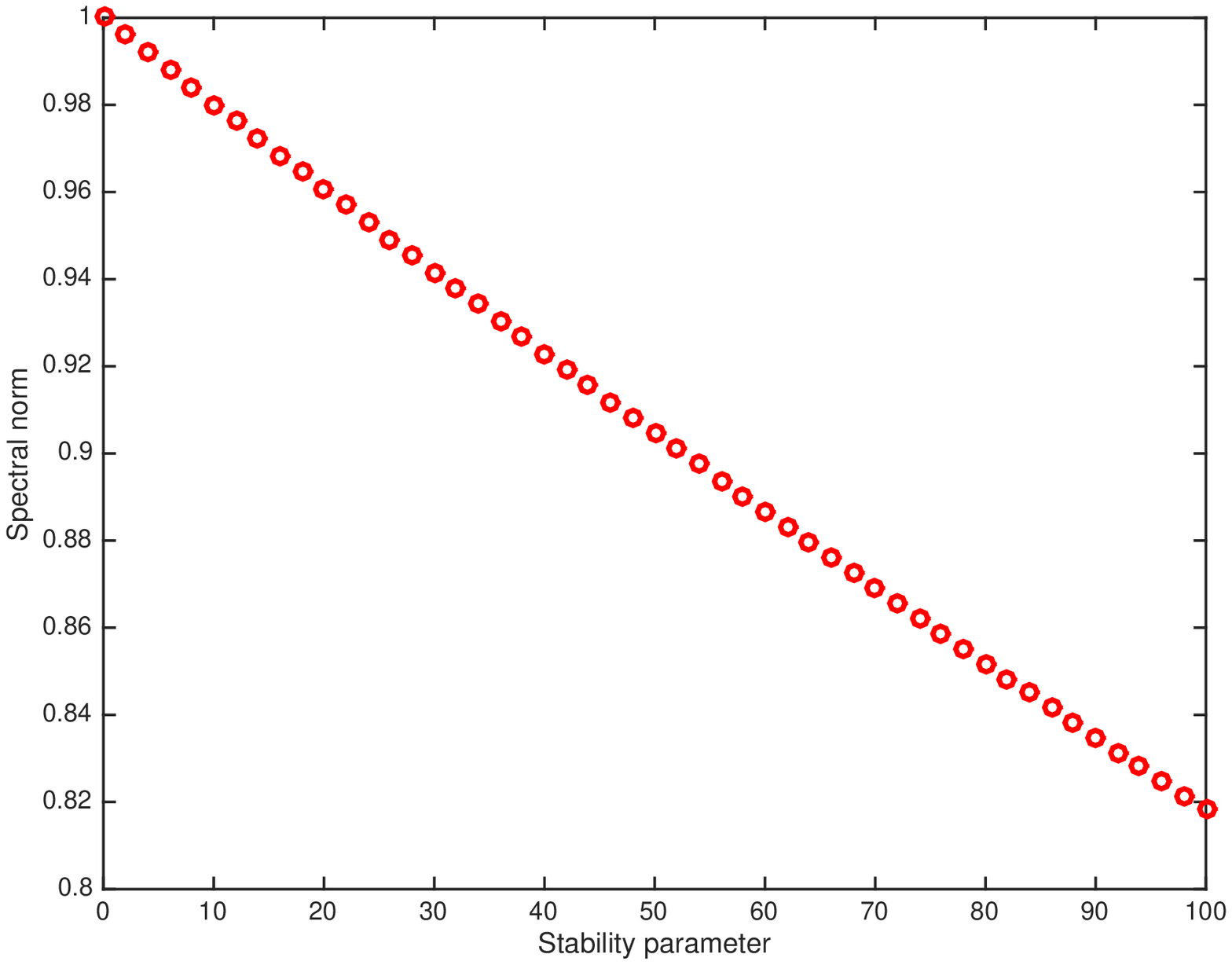}b)
\caption{Spectral  norm of the TR-BDF2 multirate method amplification matrix for system \eqref{sys1}, a) linear interpolation,
b) cubic Hermite interpolation, for increasing values of the rescaled time step.}
\label{fig:sys1_spectral}
\end{figure}

\begin{figure}[!htbp]
\centering
\includegraphics[width=0.46\textwidth]{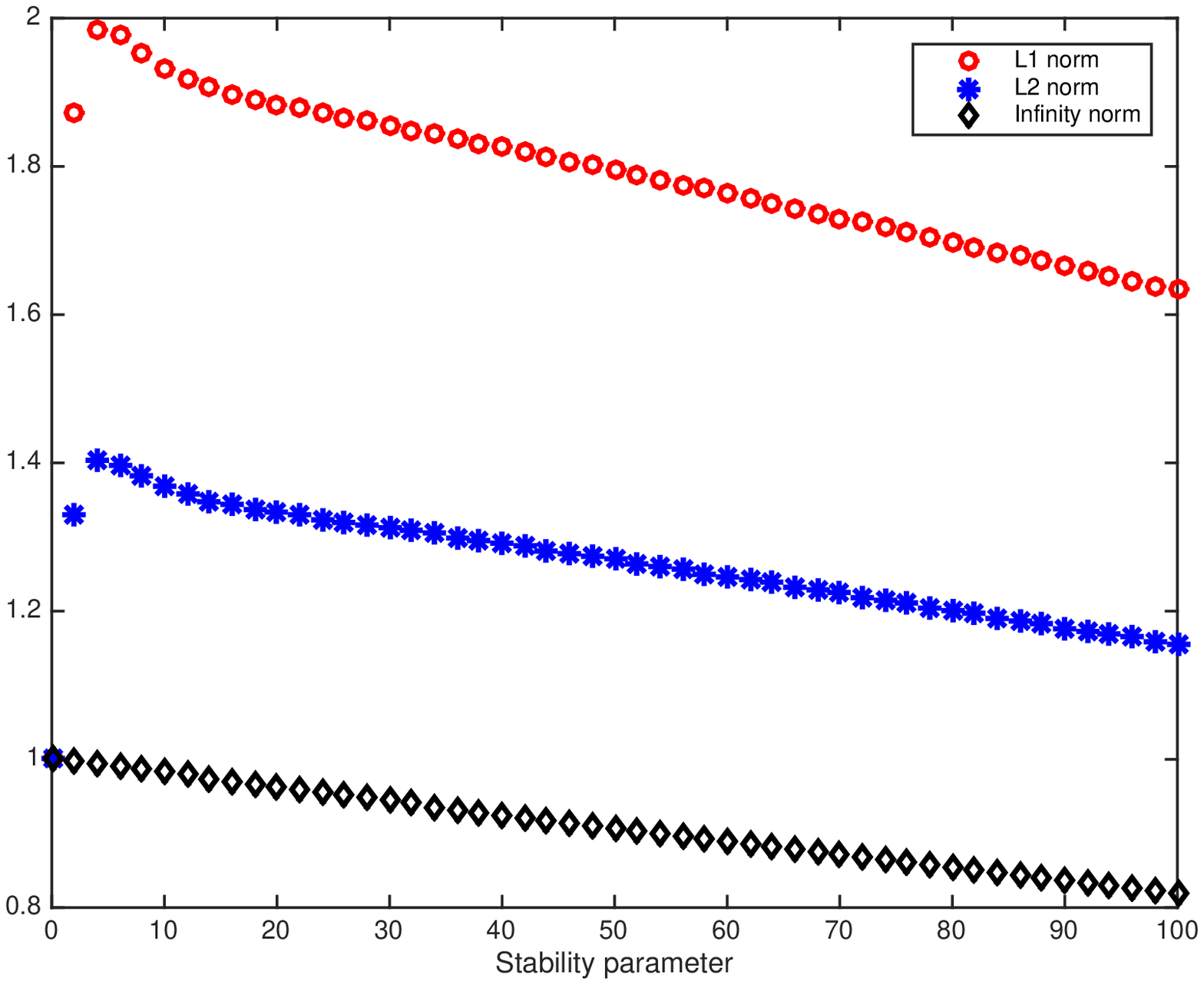}a)
\includegraphics[width=0.46\textwidth]{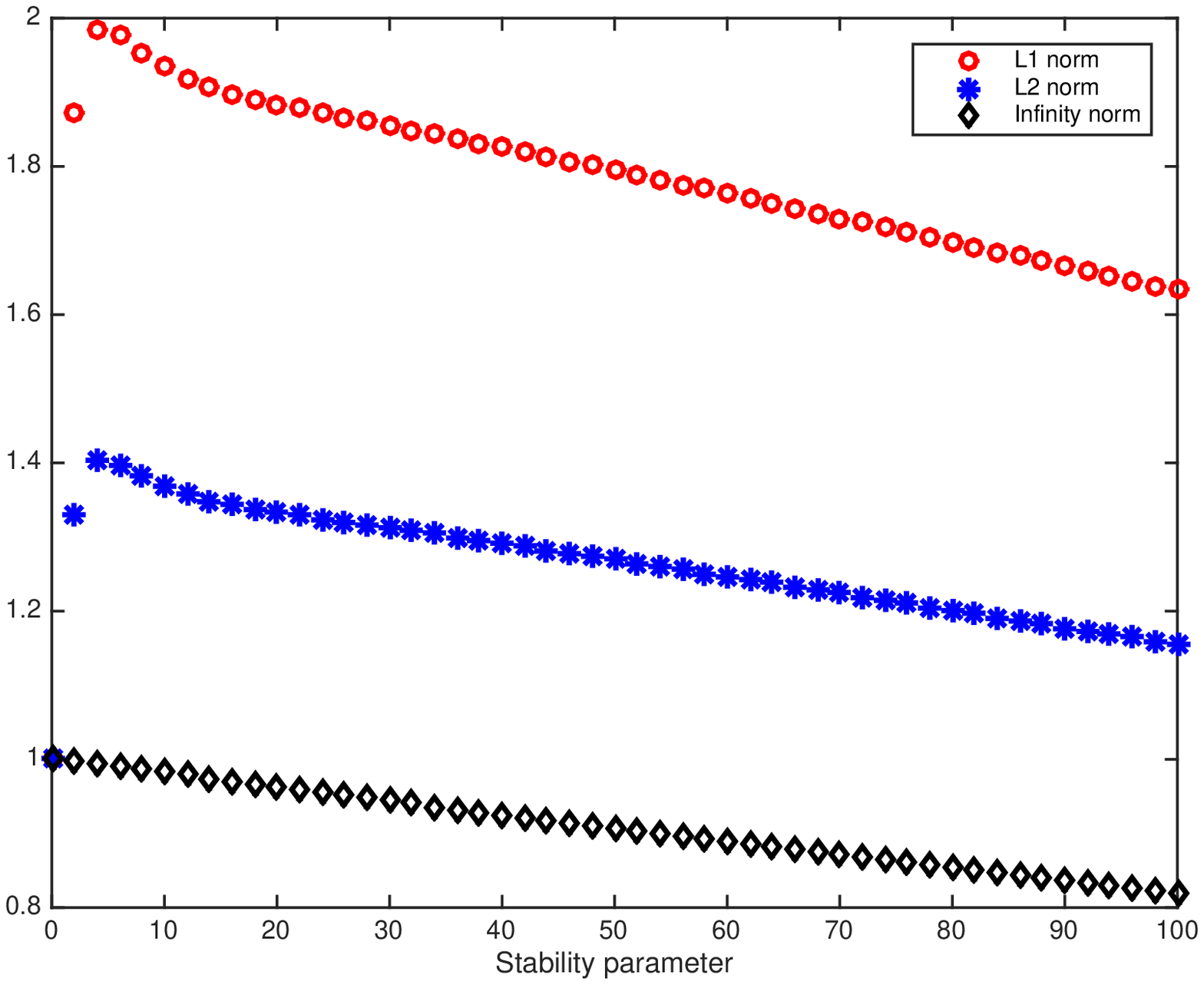}b)
\caption{Matrix norms of the TR-BDF2 multirate method amplification matrix for system \eqref{sys1}, a) linear interpolation,
b) cubic Hermite interpolation, for increasing values of the rescaled time step.}
\label{fig:sys1_norms}
\end{figure}

 Firstly, we consider a $2\times 2$ system  with the same characteristics as those 
 employed in \cite{hundsdorfer:2009}, defined by 
\begin{equation}
{\bf A}=\left[
  \begin{array}{cc}
  -1 & 1 \\
  -1000 & -1000  \\
\end{array}
  \right].
  \label{sys1}
\end{equation}
In particular, in the notation of \cite{hundsdorfer:2009}, these coefficient values correspond to the case of $\kappa=1000, $  $\beta=-1.$
The amplification matrix and its spectral and $l_1,$ $l_2,$ $l_{\infty} $ matrix norms have been computed for
different values of the time step, up to a maximum value $h_{max} $
such that $h_{max} / \max{|\lambda( \mathbf{A})|}=100.$
Therefore, the maximum time step considered is approximately 100 times larger than that of the customary
stability restriction  for explicit schemes. 
The values of the matrix norms are reported in figures \ref{fig:sys1_spectral}, \ref{fig:sys1_norms}, for the spectral
norm and the other matrix norms, respectively.
It can be observed that no stability problems arise and that the multirate algorithm appears to be insensitive to
the choice of the interpolation method.

\begin{figure}[!htbp]
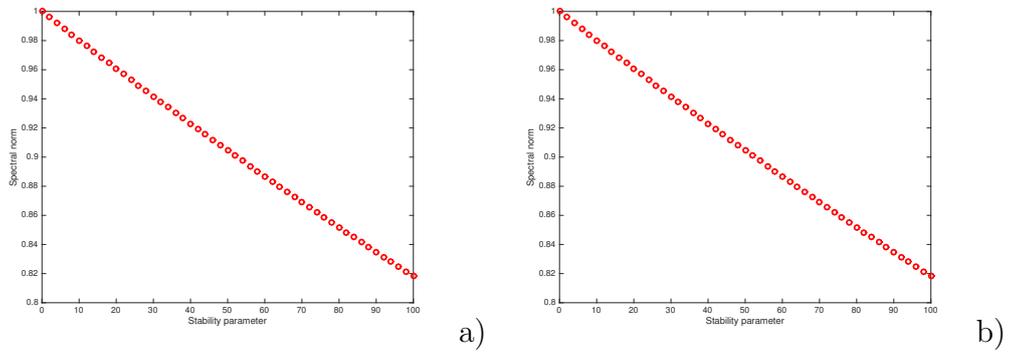

\centering
\includegraphics[width=0.46\textwidth]{./figures/savcenco_linear_spectral_trbdf2.eps}a)
\includegraphics[width=0.46\textwidth]{./figures/savcenco_cubic_spectral_trbdf2.eps}b)
\caption{Spectral  norm of the TR-BDF2 multirate method amplification matrix for system \eqref{sys2}, a) linear interpolation,
b) cubic Hermite interpolation, for increasing values of the rescaled time step.}
\label{fig:sys2_spectral}
\end{figure}

\begin{figure}[!htbp]
\centering
\includegraphics[width=0.46\textwidth]{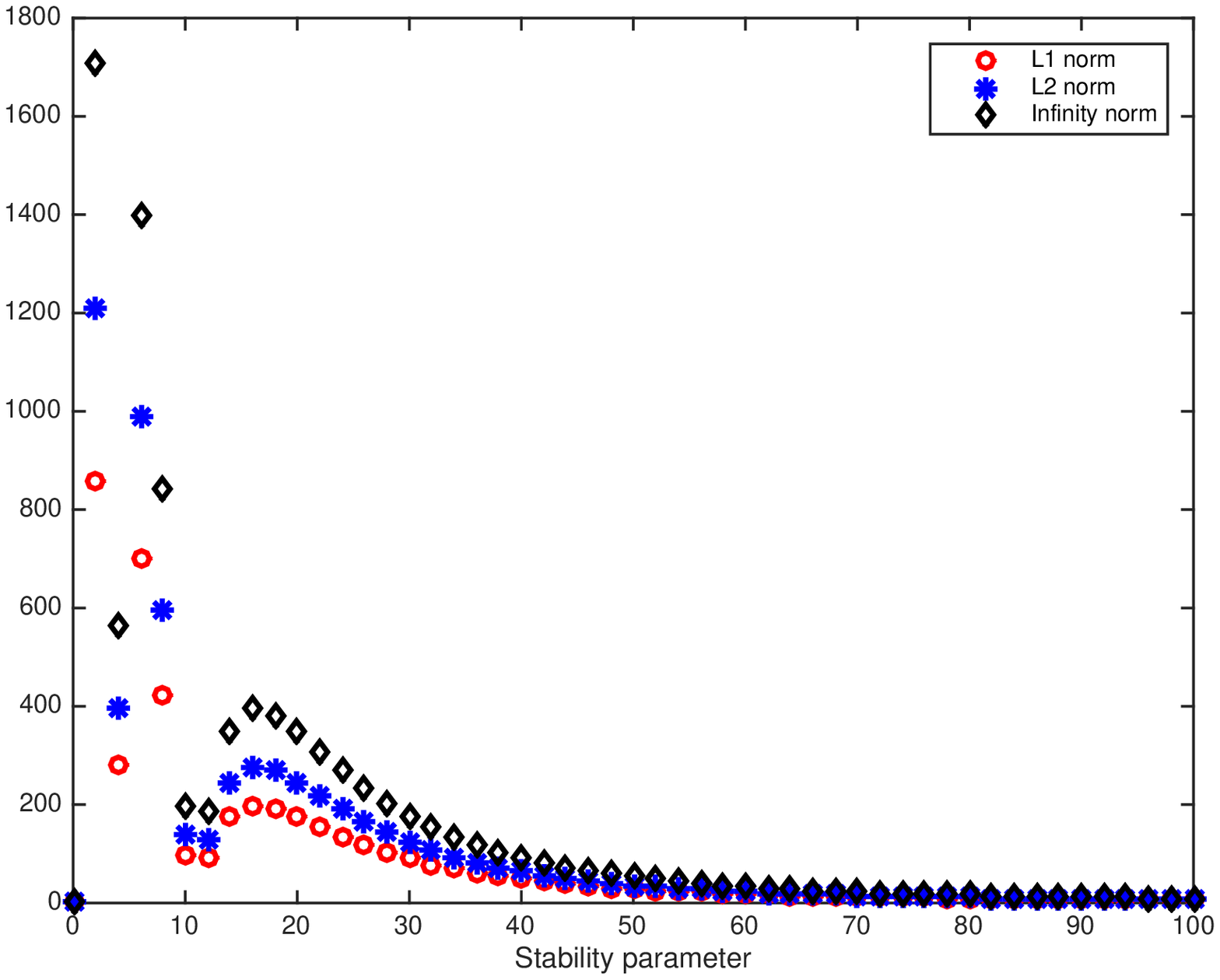}a)
\includegraphics[width=0.46\textwidth]{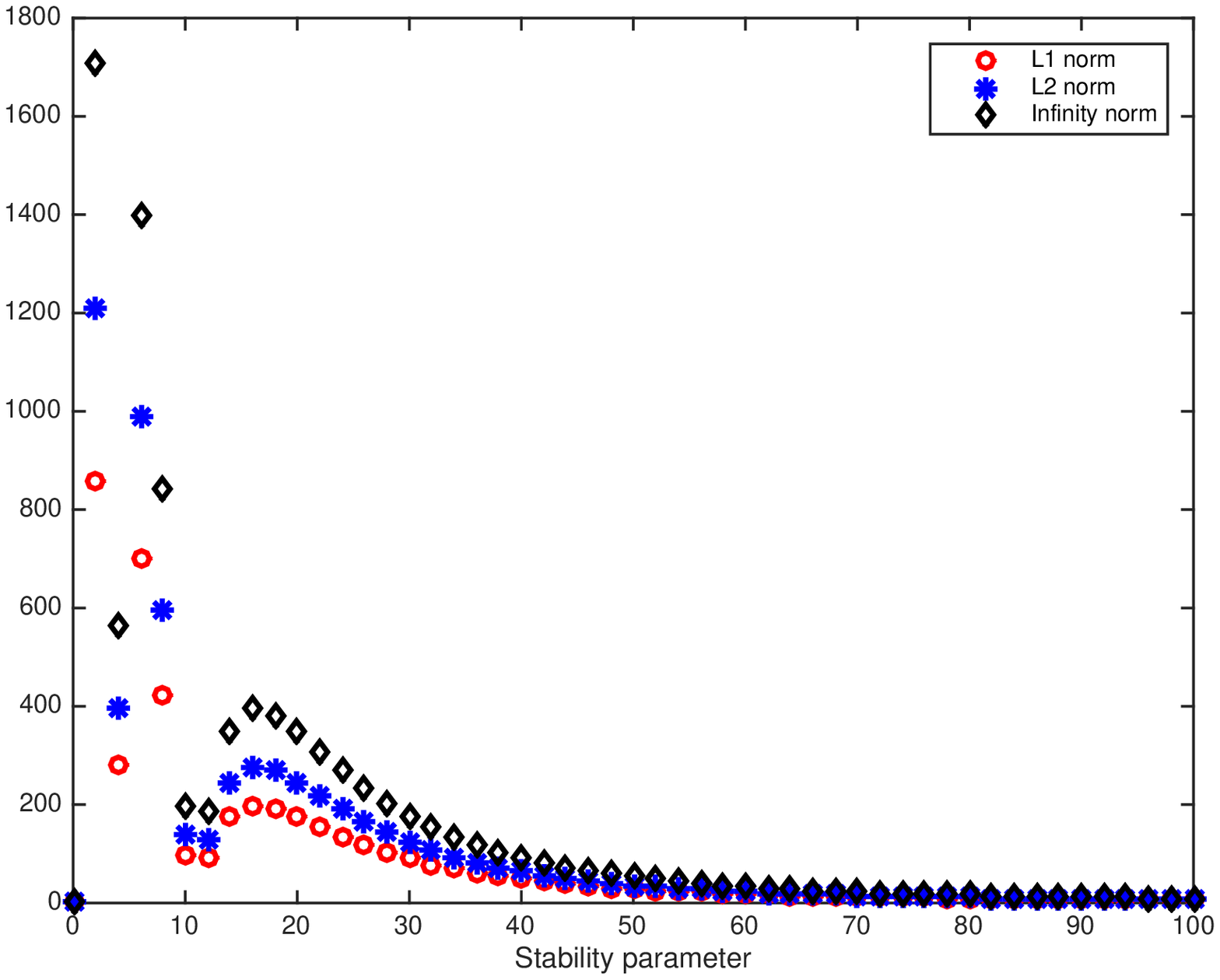}b)
\caption{Matrix norms of the TR-BDF2  multirate method amplification matrix for system \eqref{sys2}, a) linear interpolation,
b) cubic Hermite interpolation, for increasing values of the rescaled time step.}
\label{fig:sys2_norms}
\end{figure}

\begin{figure}[!htbp]
\centering
\includegraphics[width=0.46\textwidth]{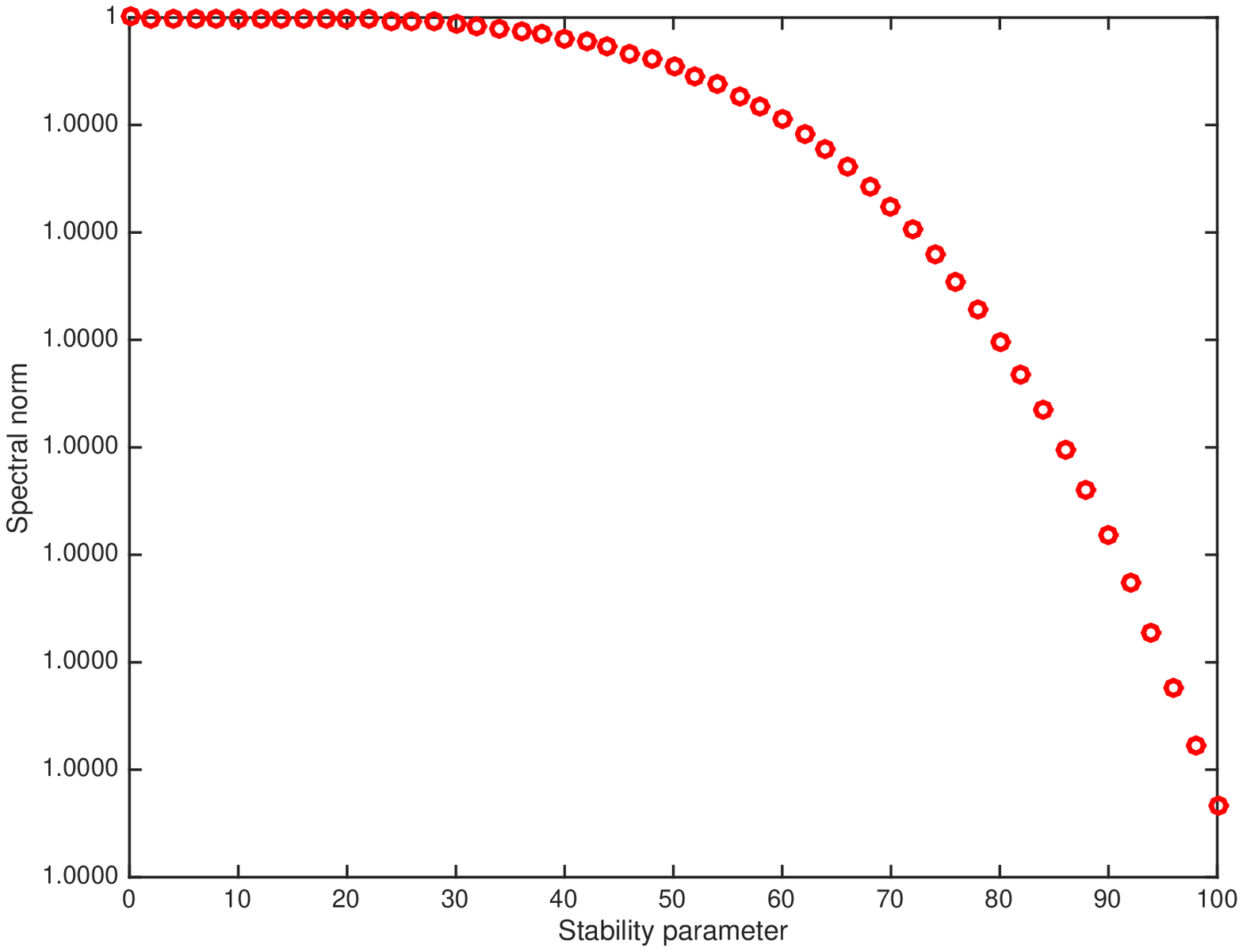}a)
\includegraphics[width=0.46\textwidth]{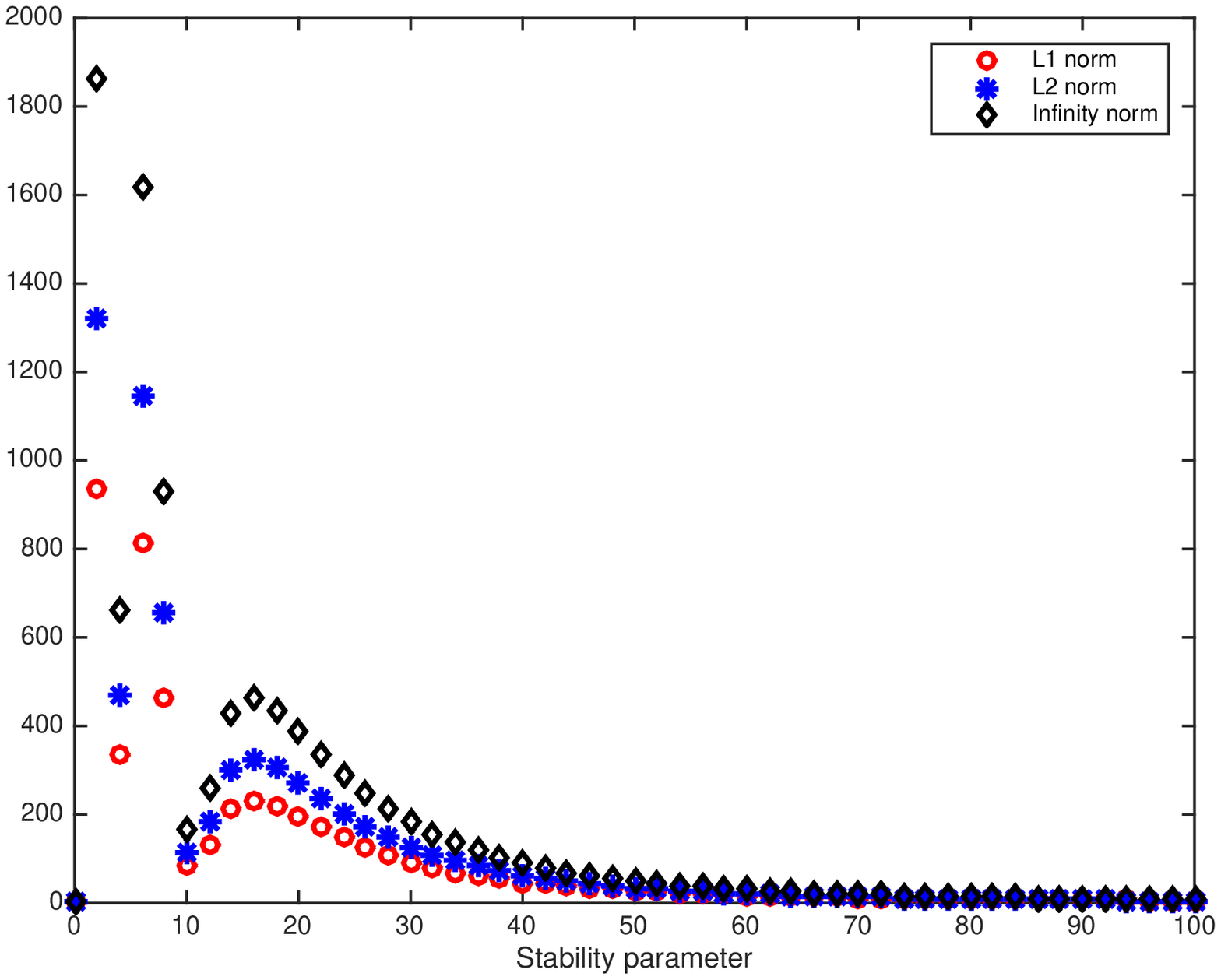}b)
\caption{Norms of the TR-BDF2  multirate method amplification matrix for system \eqref{sys2} without friction,
cubic Hermite interpolation a) spectral norm, b) other norms for increasing values of the rescaled time step.}
\label{fig:sys2_nofric}
\end{figure}

\begin{figure}[!htbp]
\centering
\includegraphics[width=0.46\textwidth]{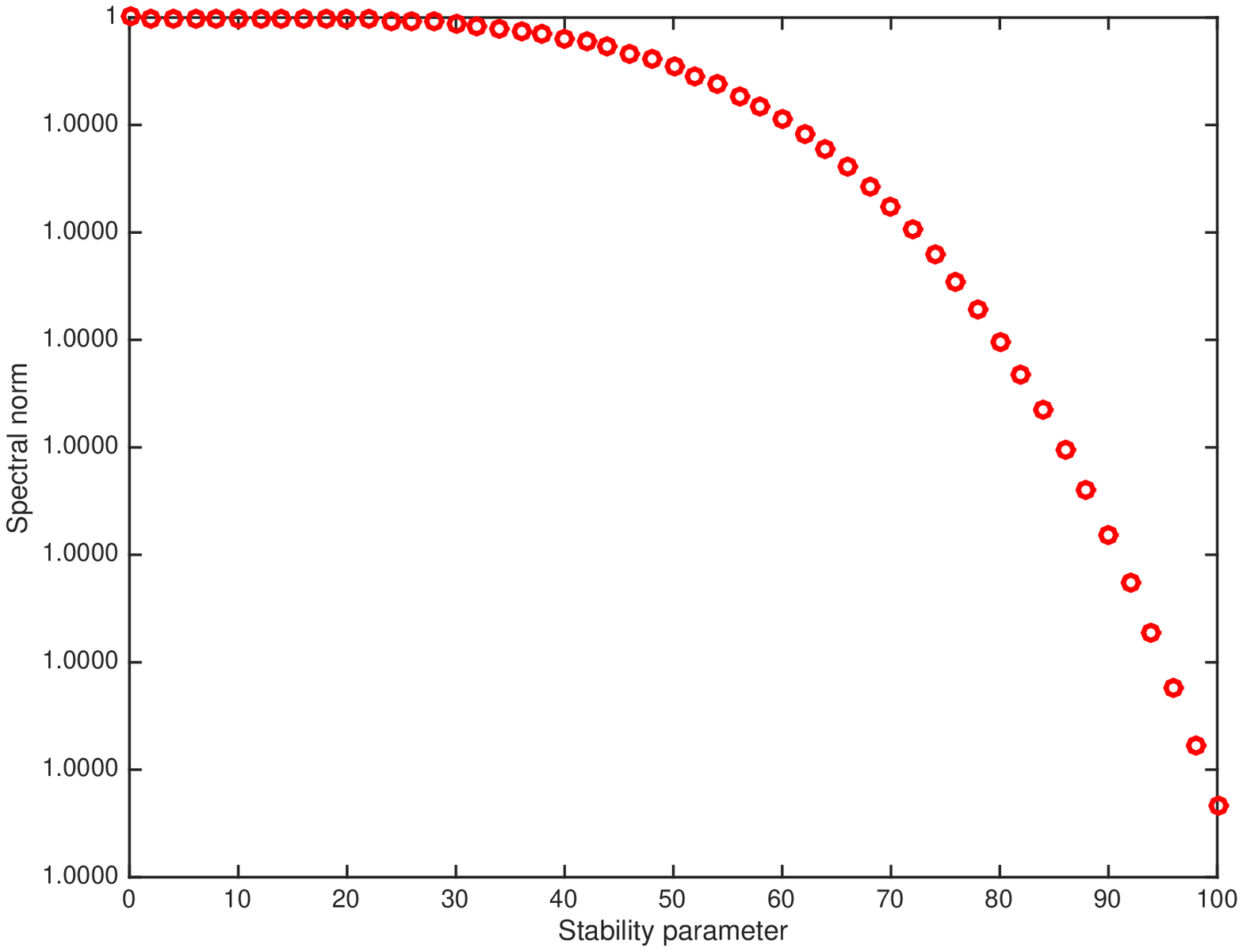}a)
\includegraphics[width=0.46\textwidth]{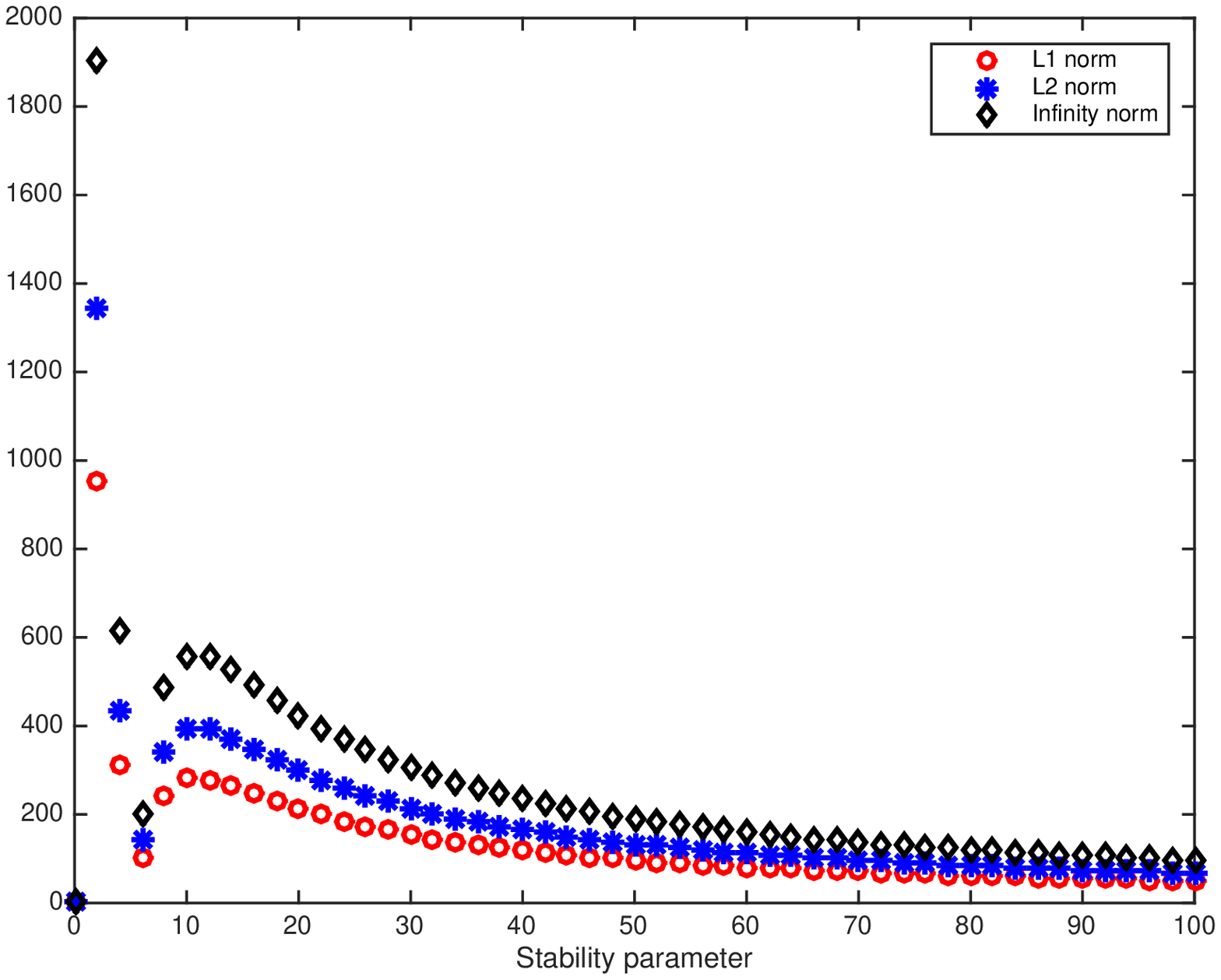}b)
\caption{Norms of the TR-BDF2  single rate method amplification matrix for system \eqref{sys2} without friction,
 a) spectral norm, b) other norms for increasing values of the rescaled time step.}
\label{fig:sys2_nofric_sr}
\end{figure}

  We have then considered the $4\times 4$ dynamical system defined
\begin{equation}
{\bf A}=\left[
  \begin{array}{cccc}
 0 & 1& 0& 0   \\
 -\frac{k_1}{m_1}& -\gamma_1 &0 &0   \\
  0 & 0& 0& 1   \\
 \frac{k_2}{m_2}&0 &-\frac{k_2}{m_2}  &-\gamma_2   
\end{array}
  \right].
  \label{sys2}
\end{equation}
  The system represents two masses $m_1,m_2 $ such that the first is tied to a wall
  by a spring of elastic constant $k_1$ and the second is tied to the first by a spring of elastic constant $k_2, $
  while $\gamma_1, \gamma_2 $ represent friction parameters. Similar systems
  have been often used to analyze empirically the stability of  numerical methods for structural mechanics,
  see e.g. \cite{erlicher:2002}. We have considered as an example the case
 $m_1=m_2=1, $ $k_1=1, $ $k_2=10^6, $ $\gamma_1=0, $ $\gamma_2=100.$
 The values of the matrix norms are reported in figures \ref{fig:sys2_spectral}, \ref{fig:sys2_norms}, for the spectral
norm and the other matrix norms, respectively.
We observe again that, even though the stability function norm can achieve large
values, the method is still stable in the spectral norm sense. 
Furthermore, it does not appear   to
be sensitive to the choice of the interpolation operator. In figure \ref{fig:sys2_nofric}, the same quantities
are plotted for the same system in the case $\gamma_2=0.$ The good stability behaviour displayed in these 
tests seems to justify further work on the application
 of this multirate approach to structural mechanics problems with multiple time scales, such as those
 considered e.g. in \cite{erlicher:2002}. As a comparison,
 we also plot in figure  \ref{fig:sys2_nofric_sr} the corresponding norms for the single rate TR-BDF2 algorithm,
 that display an entirely analogous behaviour.

\begin{figure}[!htbp]
\centering
\includegraphics[width=0.46\textwidth]{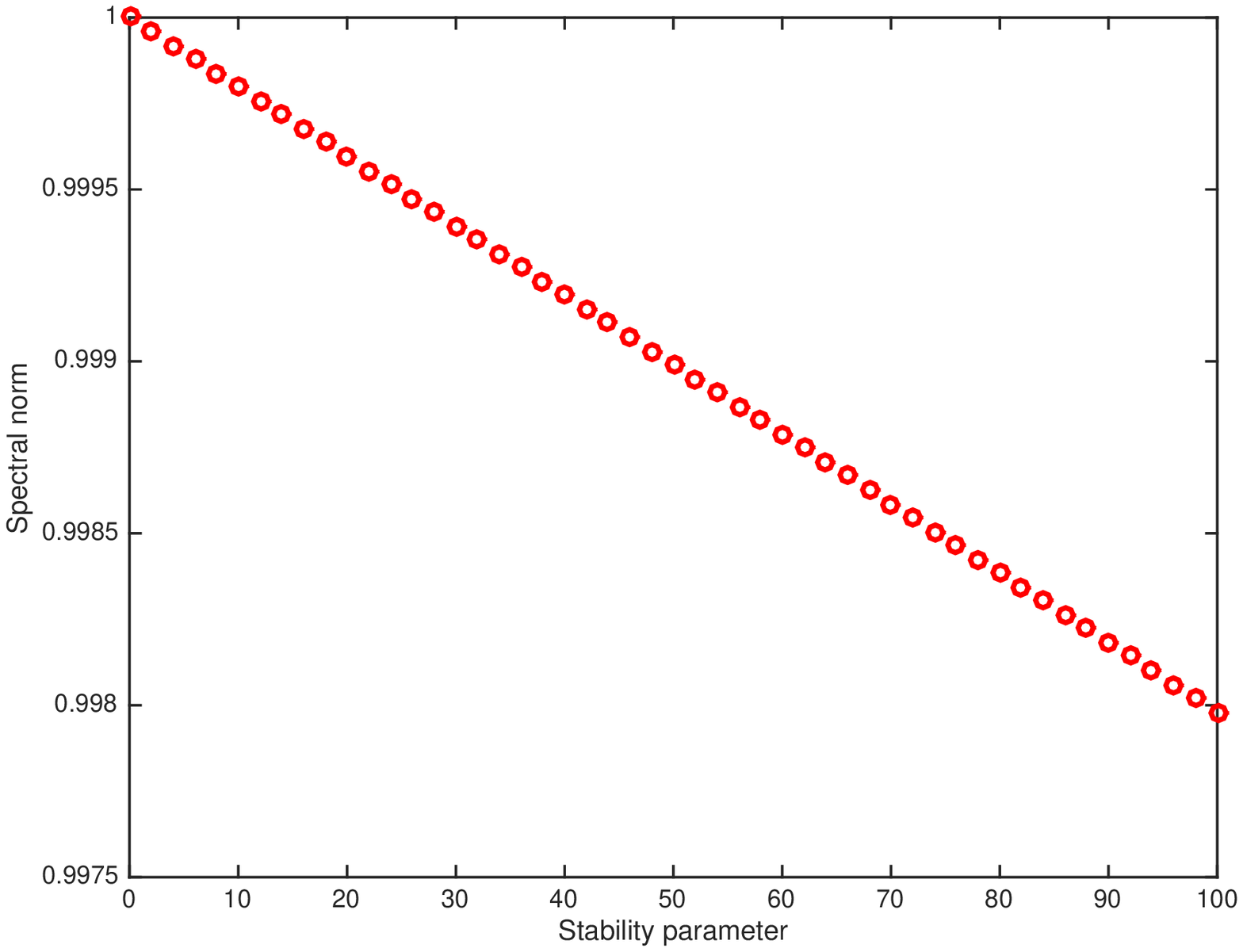}a)
\includegraphics[width=0.46\textwidth]{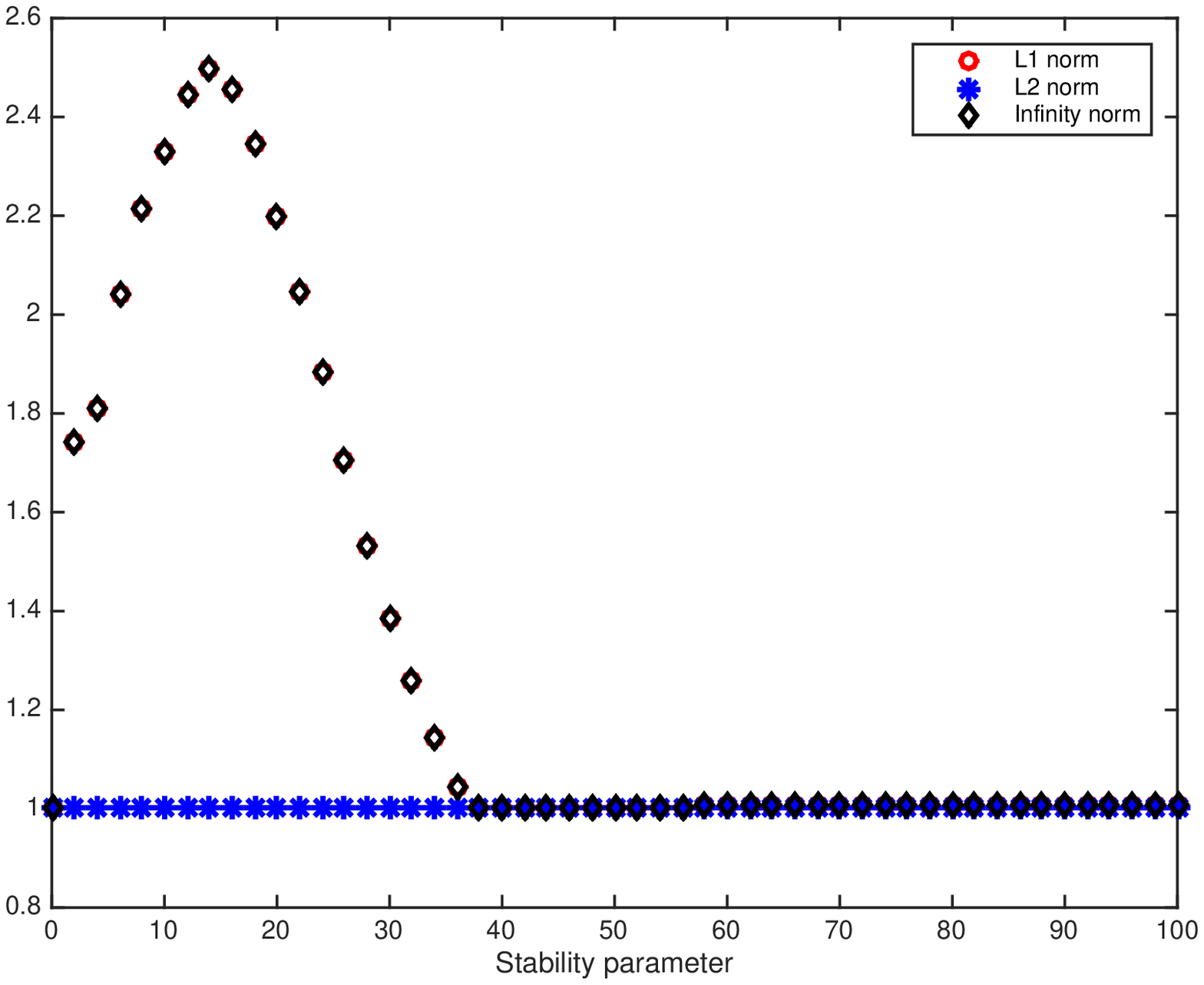}b)
\caption{Matrix norms of the TR-BDF2 multirate method amplification matrix for an 
 advection diffusion system with strongly varying
coefficients and cubic interpolation, a) spectral norm,
b) other matrix norms.}
\label{fig:sys3}
\end{figure}

\begin{figure}[!htbp]
\centering
\includegraphics[width=0.46\textwidth]{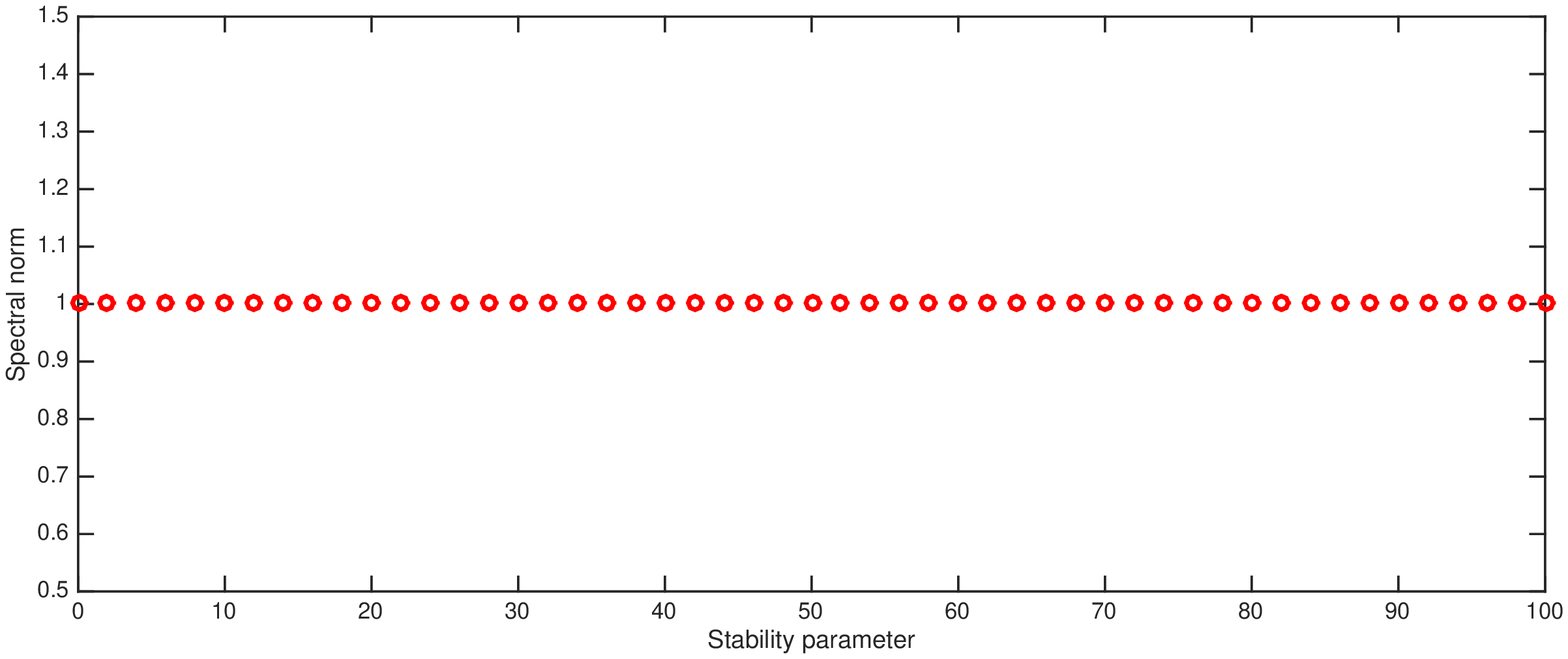}a)
\includegraphics[width=0.46\textwidth]{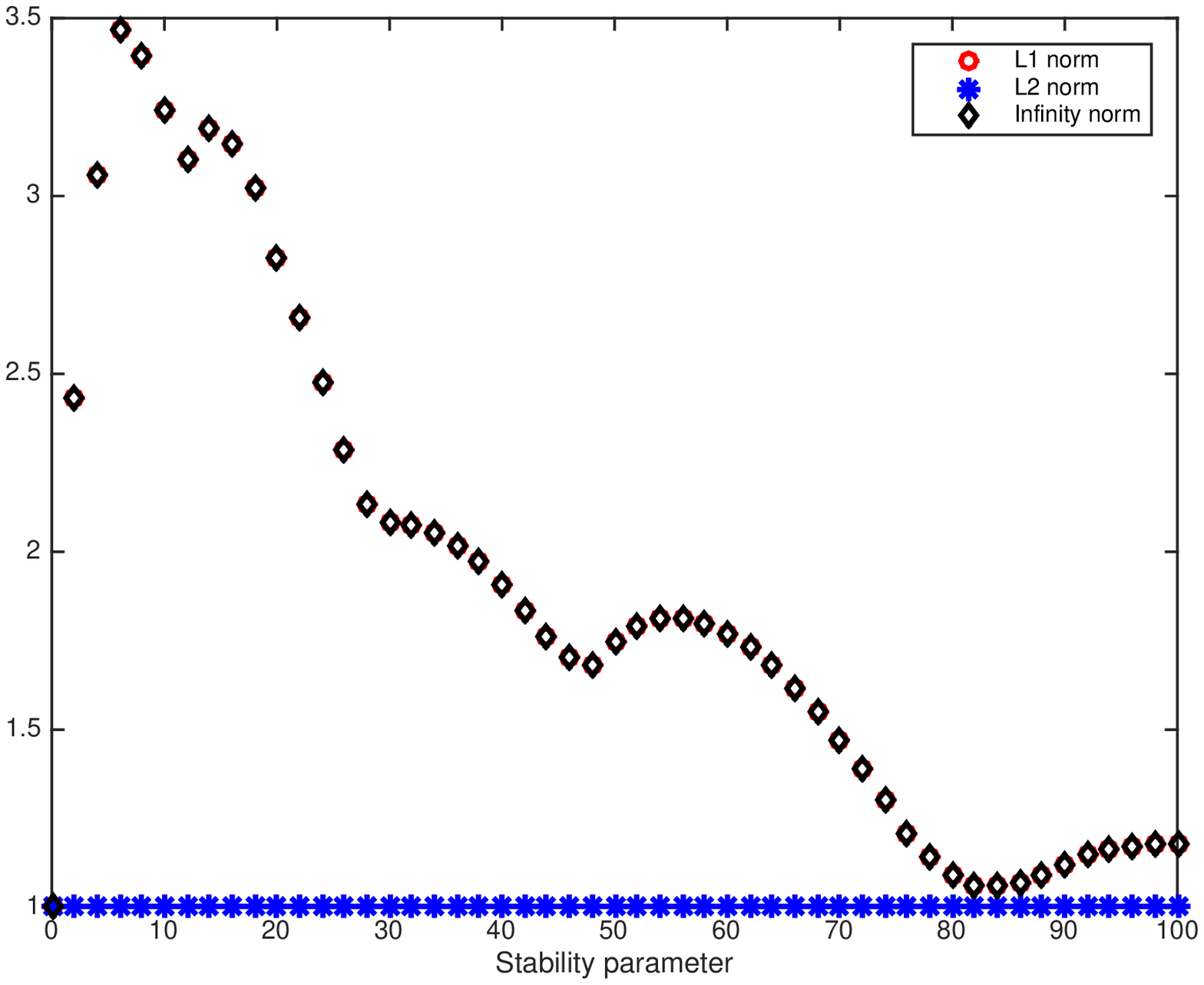}b)
\caption{Matrix norms of the TR-BDF2 multirate method amplification matrix for an 
 advection   system with strongly varying
coefficients and cubic interpolation, a) spectral norm,
b) other matrix norms.}
\label{fig:sys4}
\end{figure}

 Finally, we have considered a $40\times 40 $ linear ODE system resulting from
 the discretization of a linear heat equation by second finite differences,
 with a diffusivity parameter $10^6 $ times larger for the last 20 variables than
 for the first 20. All matrix norms are essentially identical to one over the whole range
 of time step values. Analogous results are obtained considering the corresponding
 system for linear advection-diffusion equation, taking in this case centered finite differences
 for the advection discretization and assuming the velocity and diffusivity parameters
 to be $10^4 $ times larger for the last 20 variables than
 for the first 20.  The values of the matrix norms are reported in figure \ref{fig:sys3} for the cubic interpolation case.
 In figure \ref{fig:sys4}, the analogous quantities are plotted instead in the case of a pure advection system
 with the same velocity values as in the previous case.
 The good stability properties that result in this case seem to justify further work on the application
 of this multirate approach to hyperbolic equations with multiple time scales, such as those
 considered e.g. in \cite{tumolo:2015}.

  \section{Numerical experiments}
 \label{tests} 
 \indent
 
 Several benchmark problems already considered in the literature have been employed to assess the performance of
 the multirate algorithm outlined in section \ref{selfa}.  The self adjusting multirate approach based on the TR-BDF2 solver has been
 compared to its single rate counterpart, considering several performance indicators, such as
 the CPU time required by the codes,   the number of scalar function calls and
 the number of components being computed at any given time step,  as an estimate of the computational workload for each time step.
 In all the numerical experiments, the TR-BDF2 multirate algorithm employed the cubic Hermite interpolant and a nonlinear solver
 based on the Newton method with approximate computation of the Jacobian to solve the nonlinear systems involved.
	
\subsection{The Inverter Chain Problem}
\label{sec:inverterchain}
For the first test we consider the inverter chain problem which is an important test problem from the field of electrical circuits. It has also been considered e.g.  in \cite{savcenco:2007}, \cite{verhoeven:2007}.  
The system of equations is given by

\begin{equation}
\begin{aligned}
	y_1'(t)&=U_{op} - y_1(t) -\Gamma g(u(t),y_1(t)),\\
	y_j'(t)&=U_{op} -y_j(t) - \Gamma g(y_{j-1}(t),y_j(t)), \hspace{2em} j=2,3,\ldots ,m,
\end{aligned}
\end{equation}
where $y_j, j=1,\ldots,m$  represent a chain of $m$ inverters,
 $U_{op}$ is the operating voltage corresponding to the logical value 1 and $\Gamma$ serves as a stiffness parameter. 
The function $g$ is defined by 

\begin{equation}
	g(y,z) = \left( \max\left( y-U_{\tau},0\right)\right)^2 - \left( \max\left( y-z-U_{\tau},0\right)\right)^2. 
\end{equation}
The   parameter values used were $m=500,$ $ \Gamma=100,$ $ U_{\tau}=1 $ and $U_{op}=5.$
The initial condition was defined as

\begin{equation}
	y_j(0)=\left\{ 
\begin{aligned}
       6.247\times 10^{-3}\hspace{2em}   \text{for $j$ even}\\
       5 \hspace{6em}   \text{for $j$ odd},
\end{aligned}
	\right.
\end{equation}
while the input signal was given by

\begin{equation}
	u(t) = \left\{
\begin{aligned}
	&t-5 \hspace{6em} &\text{ for } 5\leq t\leq10\\
	&5 \hspace{6em} &\text{ for } 10\leq t \leq15\\
	&\frac{5}{2}(17-t) \hspace{3em} &\text{ for } 15\leq t\leq17\\
	&0 \hspace{6em}  &\text{otherwise}
\end{aligned}\right.
\end{equation}
In this test case, the values of all the components vary in the range $(0,5)$, so that only the absolute tolerance was used. 
\begin{figure}[!htb]
	\includegraphics[width=0.45\textwidth]{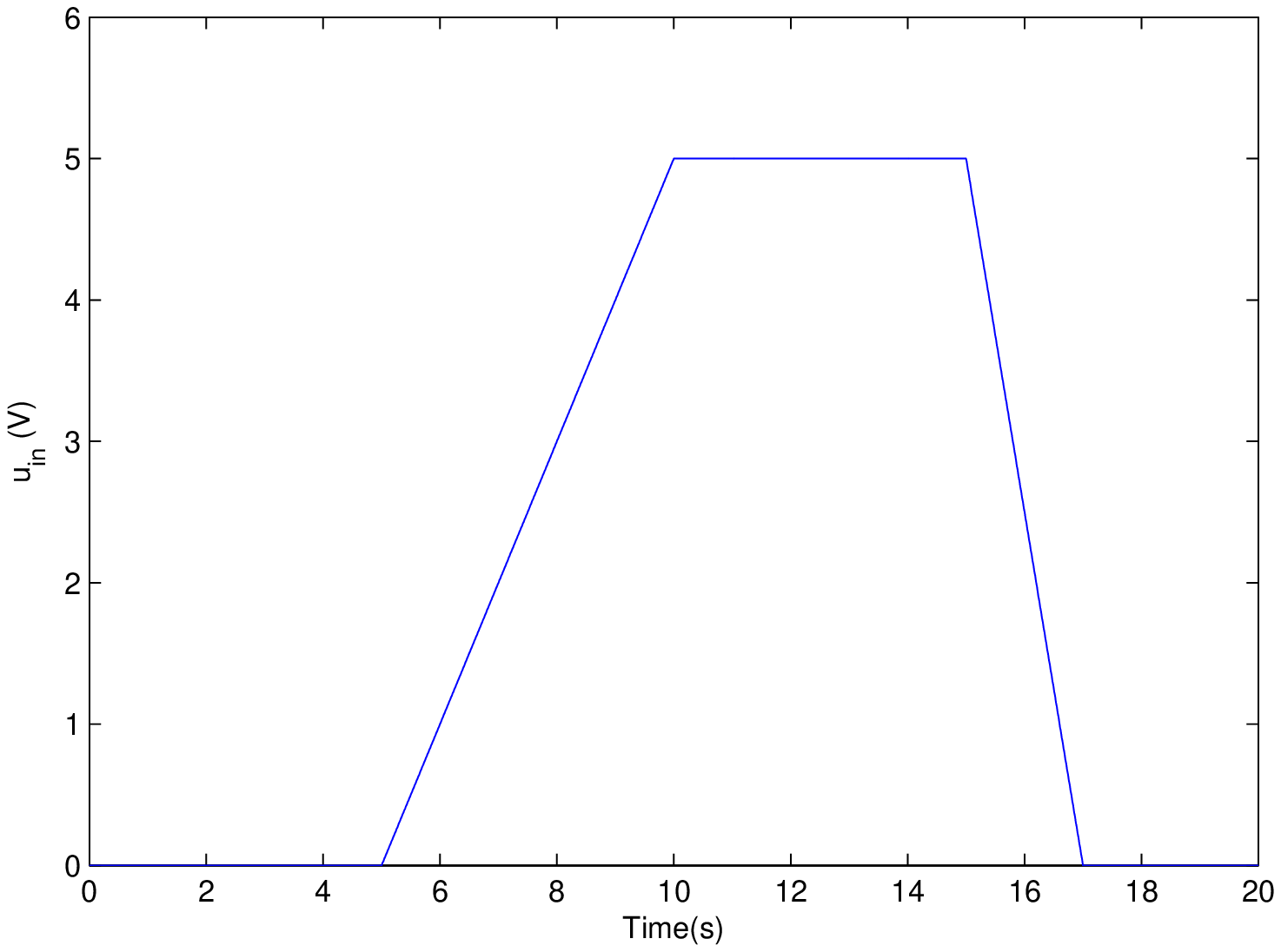}a)
	\includegraphics[width=0.45\textwidth]{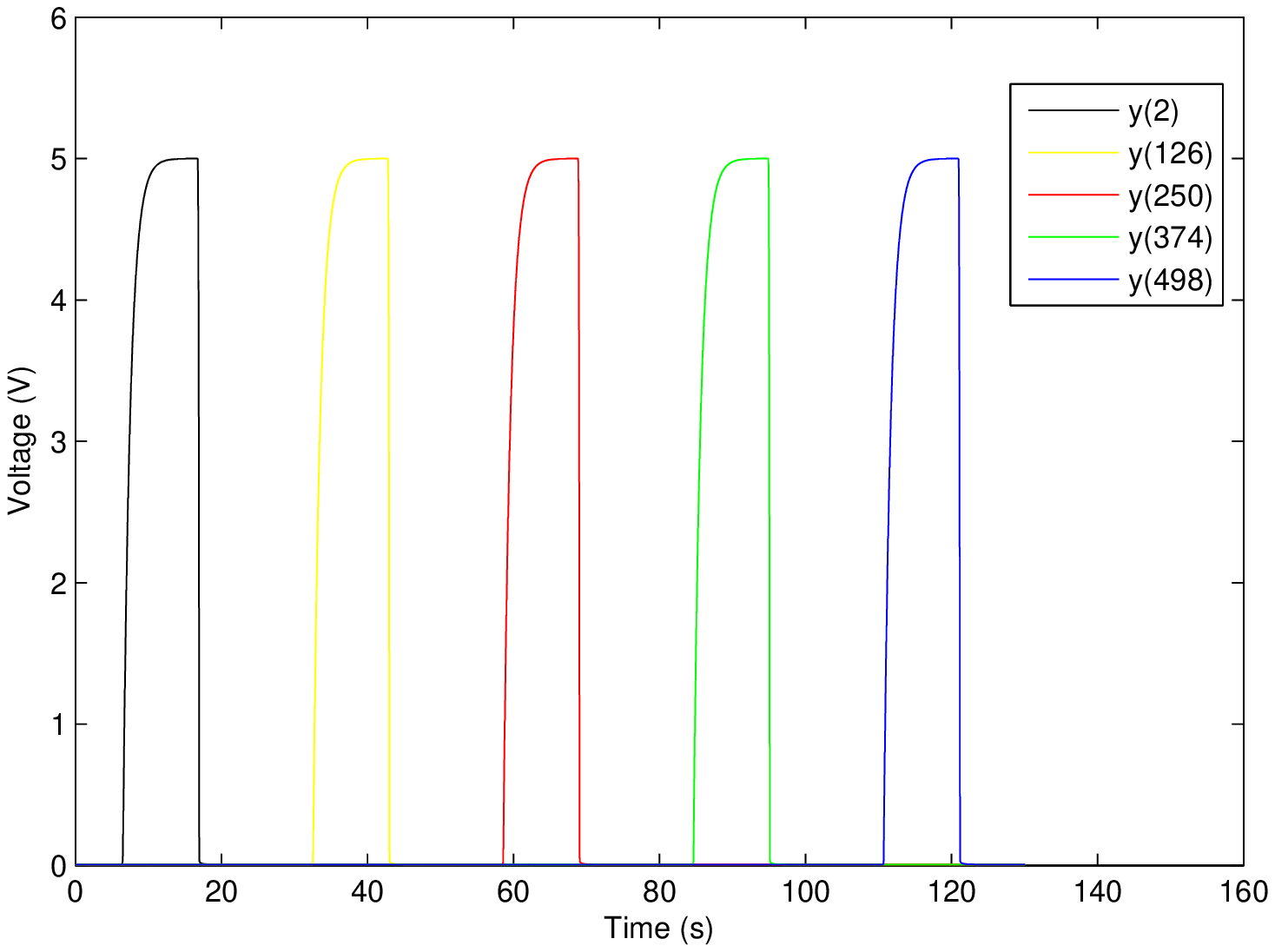}b)
	\caption{Solution of the inverter chain problem: a) input signal applied to the first inverter as a function of time,
	b) variation of voltage in 5  inverters with respect to time}
	\label{fig:INPUT_IC}
\end{figure}

\begin{figure}[!htb]
\centering
\includegraphics[width=0.7\textwidth]{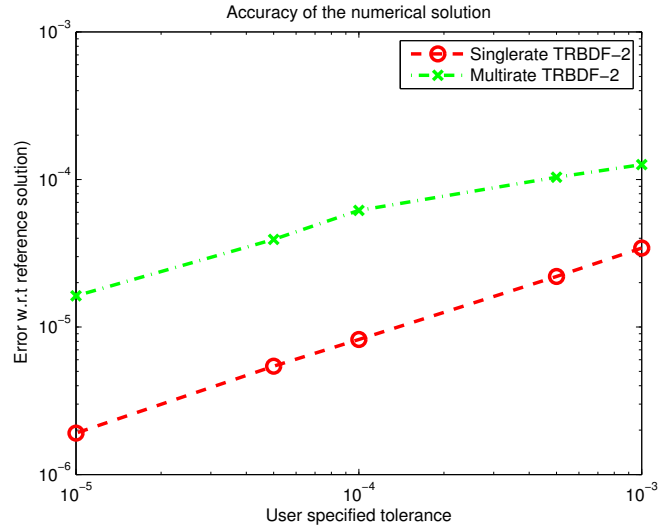}
\caption{Errors of the multirate and single rate  TR-BDF2 algorithms for the inverter chain problem, with respect to a reference solution obtained using the MATLAB \texttt{ode15s}.}
\label{fig:ACC_IC}
\end{figure}

Figure \ref{fig:INPUT_IC} shows the solution of the inverter chain problem described above.   The input signal can be seen to be propagating across the chain of inverters. It is seen that the input signal is also smoothed by the inverter chain. The reference solution was computed by the \texttt{ode15s} MATLAB solver with stringent tolerance and small maximum time step. The same reference solution
was employed to assess the accuracy of the multirate approach with respect to its single rate counterpart.  
 The maximum norm  errors at the final time of the simulation are shown in figure \ref{fig:ACC_IC}.
 Even though the   multirate  approach introduces a larger error than that of the
 single rate counterpart,  the errors committed are in general close to the specified tolerances. 

\begin{figure}[!htb]
\centering
\includegraphics[width=0.7\textwidth]{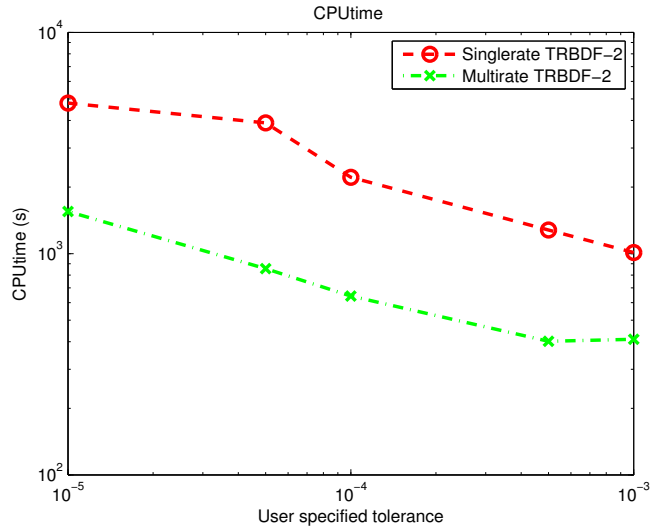}
\caption{CPU time taken by the multirate and single rate   TR-BDF2 algorithms for the inverter chain problem.}
\label{fig:CPU_HC}
\end{figure}
\begin{figure}[!htb]
\centering
\includegraphics[width=0.7\textwidth]{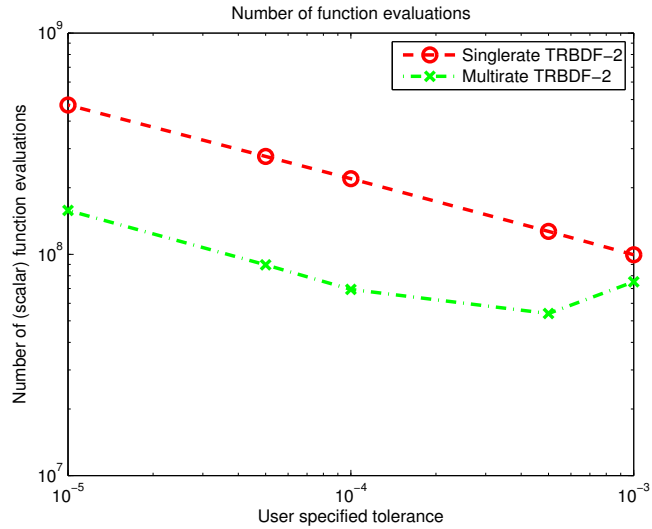}
\caption{Number of scalar function evaluations for the multirate and singlerate  TR-BDF2 algorithms for the inverter chain problem.}
\label{fig:NUM_IC}
\end{figure}

The performance of the algorithms with respect to the different metrics was analyzed for different values of the tolerance. 
Figure \ref{fig:CPU_HC} shows the plot of CPU time taken by the multirate and single rate   TR-BDF2 algorithms.  
 At a tolerance of $10^{-5}, $  a speed-up  of a factor 3 was  achieved. Figure \ref{fig:NUM_IC} shows  the total number of scalar function evaluations  at different tolerances. Again it can be seen that the performance of the multirate algorithm is better than that of
  the corresponding single rate algorithm. At a tolerance of $10^{-5}$ the number of function evaluations decreased by a factor of 3 for the TR-BDF2.   For a tolerance of $10^{-3}$, the multirate TR-BDF2 made a larger number of function evaluations than it did at tolerance $10^{-4}$. 
  Indeed, if the tolerance is not too stringent, the time step adaptation mechanism causes the algorithm to try larger steps, 
  which makes it harder for the Newton solver to converge.  

\begin{figure}[!htb]
\centering
	\includegraphics[width=0.7\textwidth]{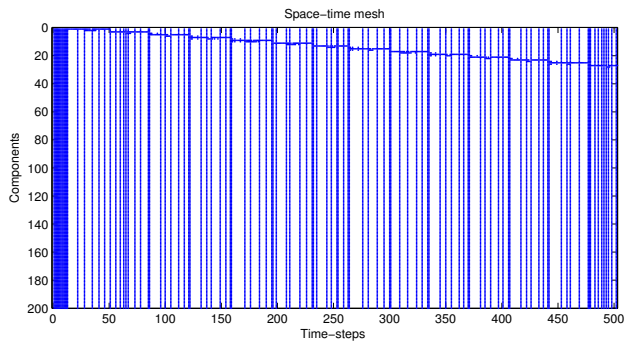}
	\caption{Space time diagram of the multirate TR-BDF2 algorithm for the inverter chain problem. Only 200 components for the first 500 time steps are shown.}  
	\label{fig:MESH_IC}
\end{figure}

\begin{figure}[!htb]
\centering
\includegraphics[width=0.7\textwidth]{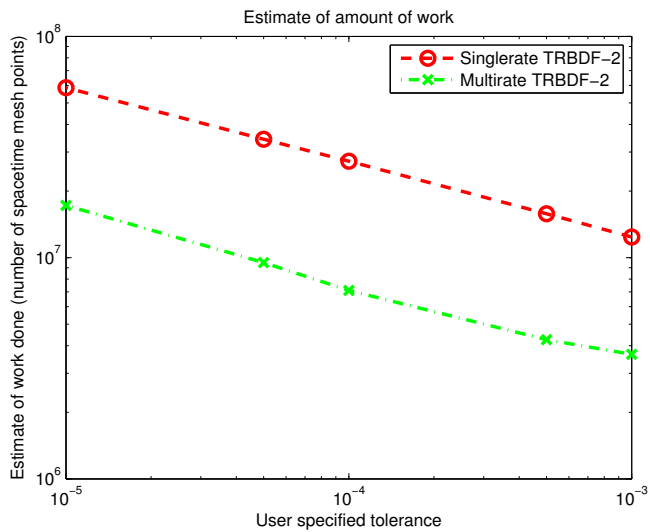}
\caption{Workload estimate for the multirate and singlerate   TR-BDF2 algorithms for the inverter chain problem.}
\label{fig:WORK_IC}
\end{figure}

In figure \ref{fig:MESH_IC}, a portion of the so-called space time diagram of the multirate algorithm is shown,
denoting all the active components at each time step. For clarity,  only a zoomed in view  of 200 inverters for the first 500 time steps
is shown. SInce the simulation is started with a very conservative value of the time step,
 the first few time steps involve all the degrees of freedom. The time step is then
 gradually increased by the   adaptation mechanism, so that at a later stage there are only a few components being refined at each time step. 
 The set of refined components moves forward across the inverters as the simulation goes on. This is as expected, because the signal is propagating across the inverter chain and the  active components   are being refined, while the others are being integrated with larger time-steps.  The number of space-time points in the   figure \ref{fig:MESH_IC} can be regarded as a measure of the computational workload of the algorithm.  Obviously, for a single rate algorithm this would simply be equal to the number of components times the number of time-steps. Figure \ref{fig:WORK_IC} shows this estimate of computational workload for the algorithms for the inverter chain problem. With regard to this metric, using a tolerance of $10^{-5}$ the multirate algorithm leads to a gain in efficiency of 3.4 times with respect to the single rate algorithm.

%

\subsection{A reaction diffusion problem}
\label{sec:rxndfsn}
\begin{figure}[!htb]
\centering
	\includegraphics[width=0.7\textwidth]{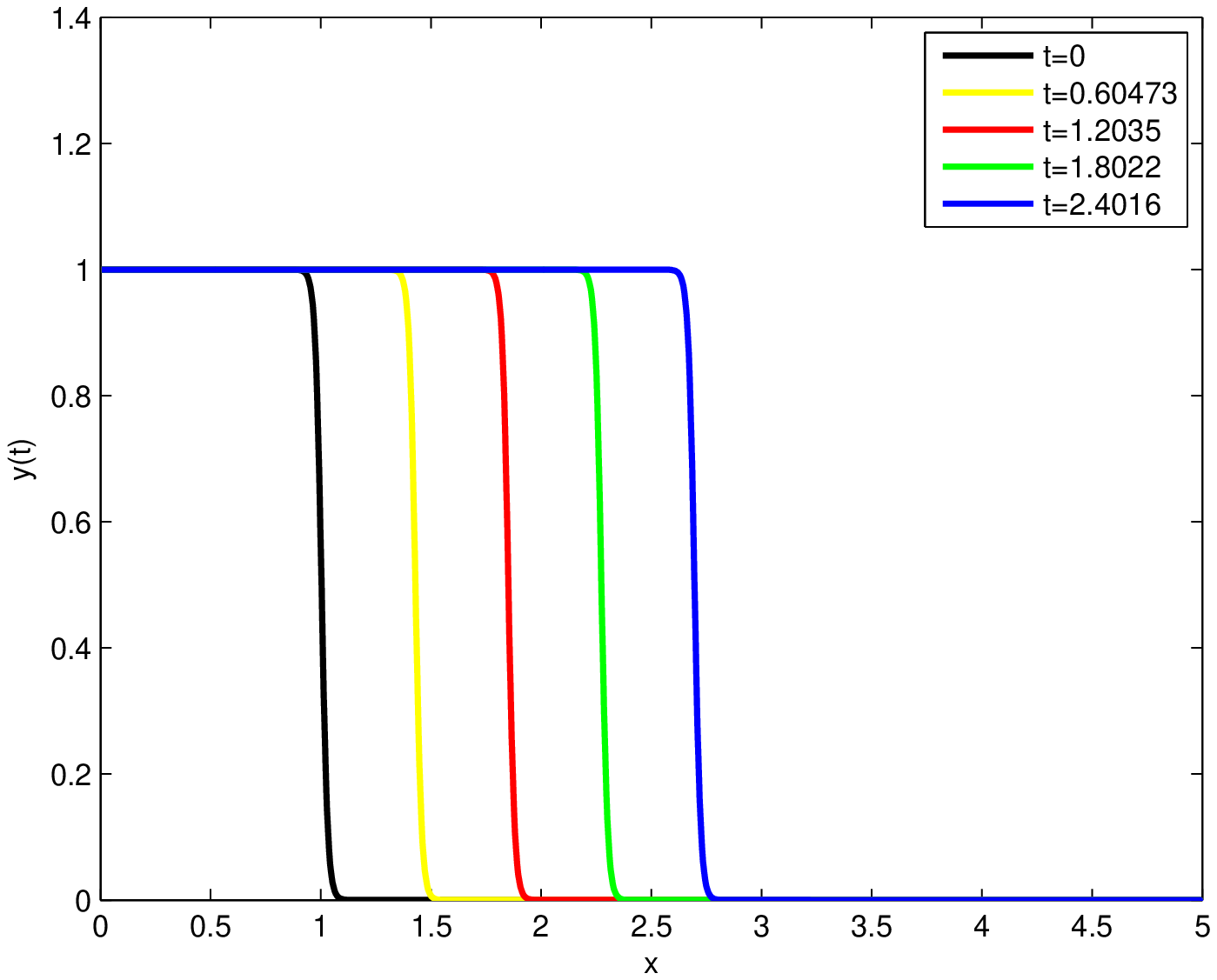}
	\caption{Reference solution of the reaction-diffusion problem  obtained by the MATLAB \texttt{ode15s} solver, shown
	 at five evenly spaced time instants in the simulation interval.}
	\label{fig:SOL_RD}
\end{figure}

We then consider a test in which the ODE problem to be solved results from the space discretization
of a partial differential equation. In particular, we consider the reaction diffusion problem also
discussed in \cite{savcenco:2007}, whose associated PDE can be written as

\begin{equation}
\label{ref:rxndfsn2}
	\frac{\partial y}{\partial t}=\epsilon\ \frac{\partial y^2}{\partial x^2} + \gamma y^2(1-y),
\end{equation}
in the domain $0<x<L$, $0<t<T$. The initial and boundary conditions are given by
$$
y(0,x)= \frac{1}{1+\exp\left(\lambda(x-1)\right)} \hspace{2em} \frac{\partial y(t,0)}{\partial x} =   \frac{\partial u(t,L)}{\partial x}  =0.
$$
The values used for parameters were $\lambda=\dfrac{1}{2}\sqrt{\dfrac{2\gamma}{\epsilon}}$, $\gamma=\dfrac{1}{\epsilon}=100$, $L=5$ and $T=3$. 
The equation was discretized in space by simple second order finite differences on a mesh of constant spacing
$\Delta x=L/N, $ with different values of $N.$ Even though no analytic solution is available,  the solution is well known
to consist in a   wave moving in the positive direction and connecting the two stable states
of the nonlinear reaction potential. Figure \ref{fig:SOL_RD} shows the solution as a function of space at different time instants.

\pagebreak

\begin{figure}[!htb]
\centering
	\includegraphics[width=0.7\textwidth]{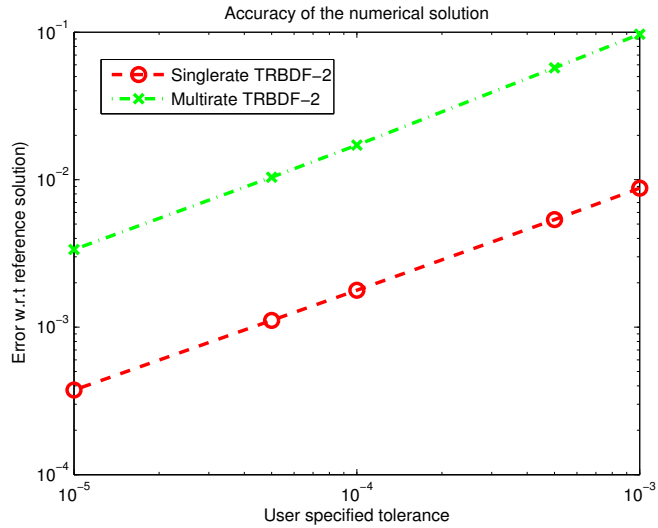}
	\caption{Errors committed by the multirate and single rate  TR-BDF2 algorithms for different tolerances for the reaction diffusion problem described in section \ref{sec:rxndfsn}}
	\label{fig:ACC_RD}
\end{figure}

\pagebreak

\begin{figure}[!htb]
\centering
	\includegraphics[width=0.45\textwidth]{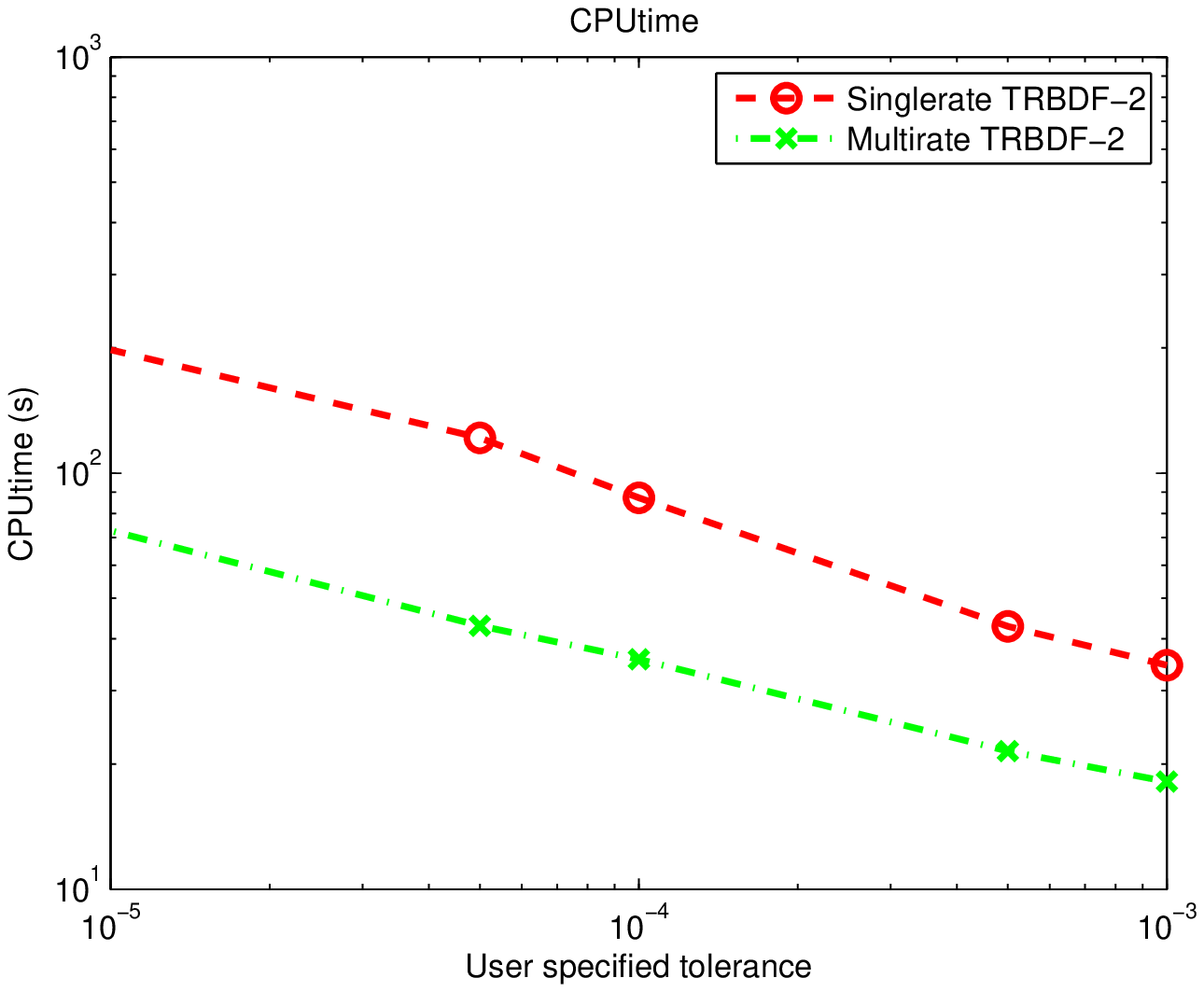}a)
	\includegraphics[width=0.45\textwidth]{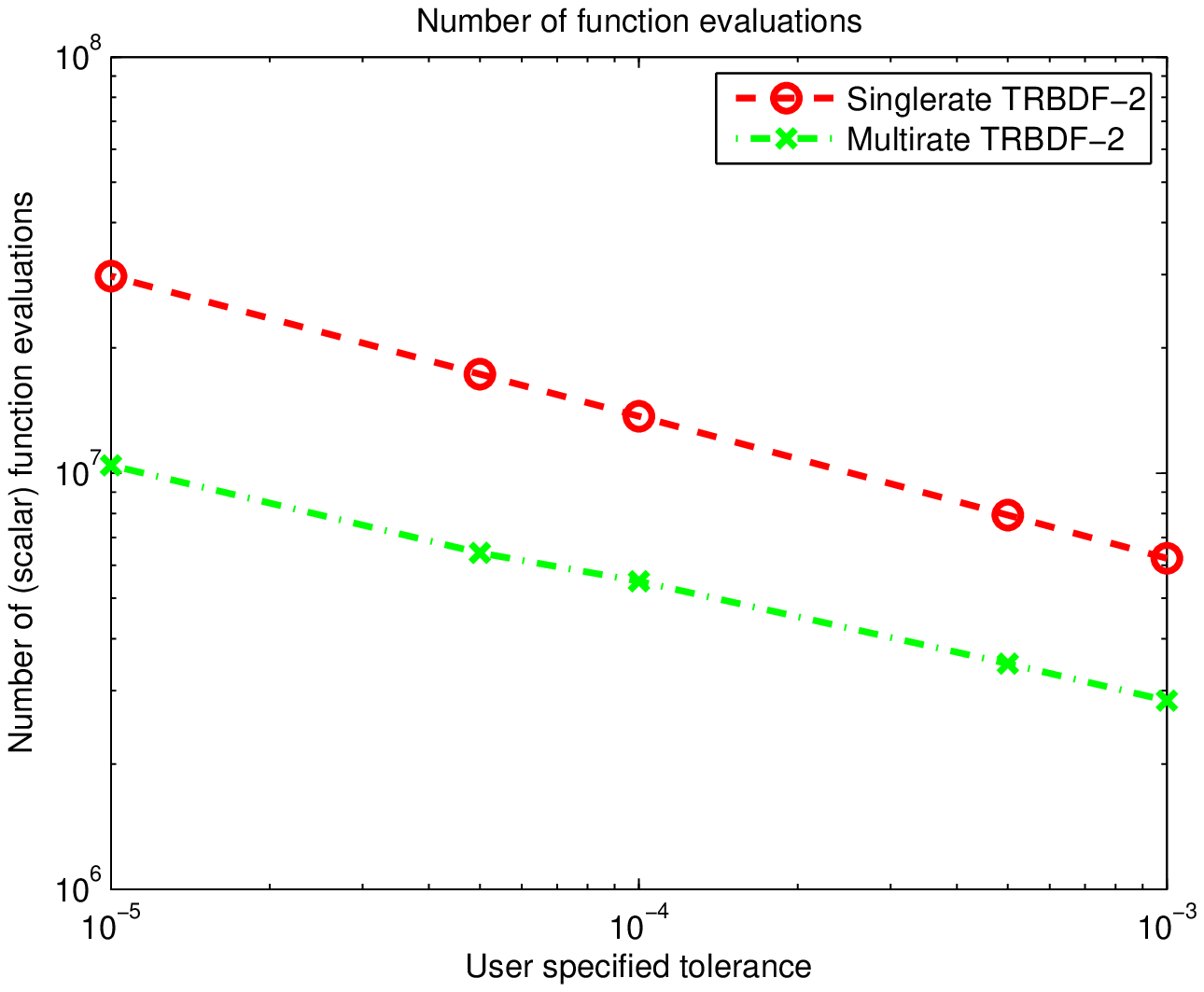}b)
	\caption{a) CPU time  and b) number of scalar function evaluations for the multirate and single rate   TR-BDF2 algorithms for different tolerances for the reaction diffusion problem.}
	\label{fig:CPU_RD}
\end{figure}

As in the previous section, we have compared the accuracy of the proposed algorithm with respect to the reference solution computed by the MATLAB \texttt{ode15s} solver with tight tolerance  and small maximum time step values.
Figure \ref{fig:ACC_RD} shows the plot of the errors at the final simulation time $T=3 \rm s $ for different values of user specified tolerance.
 The performance of the algorithms with respect to the different metrics was analyzed for different values of the tolerance. 
Figure \ref{fig:CPU_RD} shows the plot of CPU time and the number of function evaluations required by the multirate and single rate   TR-BDF2 algorithms, respectively.  
 At a tolerance of $10^{-5},$ the speed-up of the multirate algorithm with respect to the single rate algorithm was     2.7.

 \subsection{Linear Advection Equation}
 
\begin{figure}[!htb]
\centering 
\includegraphics[width=1.1\textwidth]{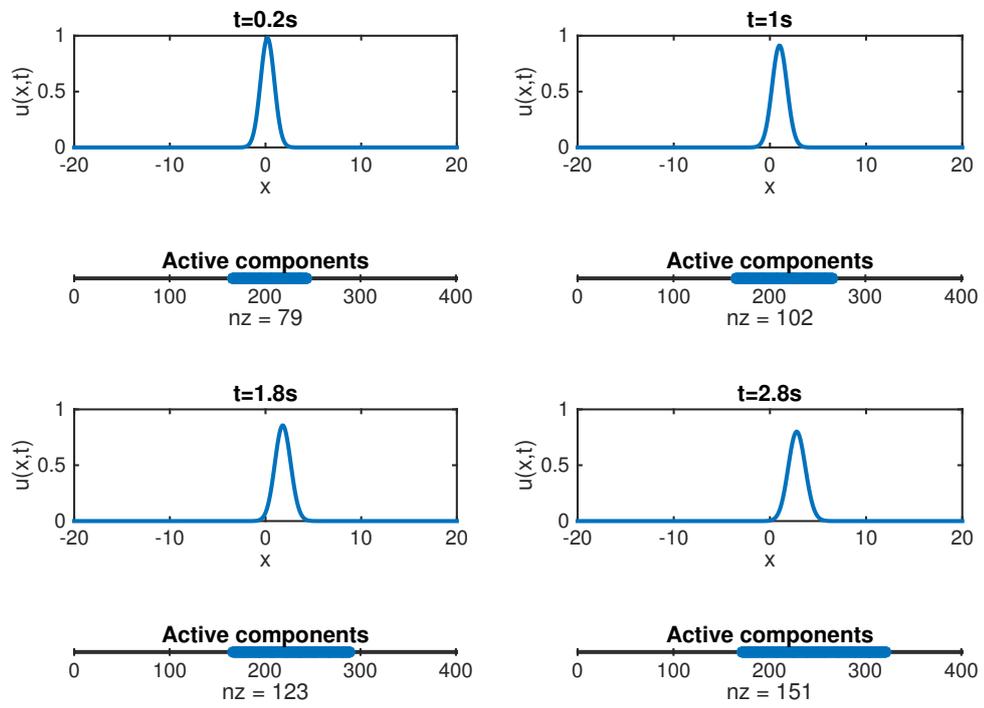}
\caption{Multirate TR-BDF2 method  with cubic interpolation for the linear advection problem.}
\label{Upwind}
\end{figure}

 As first example of hyperbolic equations we consider the linear problem

\begin{equation}
\dfrac{\partial u}{\partial t} + \dfrac{\partial u}{\partial x} = 0, \qquad x\in[-20,20]
\end{equation}
with a Gaussian initial datum.
To discretize in space we used a first order upwind method 
with periodic boundary conditions.
The interval time is $[0,3] $ with a number of cells equal to $400. $  The relative error tolerance  was
set to $ 10^{-6}, $ while the absolute error tolerance was set to
$10^{-8}. $ The initial size of the time step was taken to be $10^{-2}$.
The upwind discretization entails a relevant amount of numerical diffusion,  
which is responsible for the spreading of the initial profile. This effect also entails  that the number of computed components increases as time progresses, see Figure \ref{Upwind}.

 \pagebreak
 
 \pagebreak

%
%
%

\begin{table}[!htb]
\centering
\begin{tabular}{|c|c|c|}
\hline 
Method & CPU time  & Number of time steps \\ 
\hline 
Single rate TR-BDF2 & $21.41 s$ & $323$ \\ 
\hline 
Multirate TR-BDF2 & $11.59 s$ & $433$ \\ 
\hline 
\end{tabular} 
\caption{CPU time and number of time steps for the linear advection test.}
\label{SvsMhy}
\end{table} 
%

\begin{table}[!htb]
\centering
\begin{tabular}{|c| p{2.5cm}| p{2.5cm} | p{2.5cm}|}
\hline
Time [s] & Cubic Interp.  & Linear Interp. & Single rate \\ 
\hline
$0.2$ & $5.7\times10^{-2}$ & $5.71\times 10^{-2}$ &$5.88\times 10^{-2}$ \\ 
\hline
$1$ & $2.73\times 10^{-1}$&  $2.74\times 10^{-1}$ & $2.63\times 10^{-1}$  \\ 
\hline 
$1.8$ & $4.52\times10^{-1}$ & $4.53\times 10^{-1}$ & $4.46\times10^{-1}$ \\ 
\hline 
$2.8$ & $6.41\times10^{-1}$ & $6.43\times 10^{-1}$ & $6.41\times10^{-1}$ \\   
\hline
\end{tabular} 
\caption{Relative error in $l^{\infty}$ norm at different times for the linear advection problem. The exact solution
was used as a reference  and
the solution was computed with a) the multirate method with cubic Hermite interpolator b) the  multirate method with linear interpolator  and c) the  single-rate TR-BDF2 method.}
\label{error_linearupwind_ex}
\end{table}


\begin{table}[!htb]
\centering
\begin{tabular}{|c|p{2.5cm}| p{2.5cm}| p{2.5cm} |}
\hline
Time [s] & Cubic Interp.  & Linear Interp. & Single rate   \\ 
\hline
$0.2$ & $1.41\times10^{-6}$ & $3.07\times 10^{-5}$ &$1.38\times 10^{-6}$  \\ 
\hline
$1$ & $5.45\times 10^{-6}$&  $1.54\times 10^{-5}$ & $4.76\times 10^{-6}$  \\ 
\hline 
$1.8$ & $8.78\times10^{-6}$ & $1.49\times 10^{-5}$ & $8.80\times10^{-6}$  \\ 
\hline 
$2.8$ & $1.24\times10^{-5}$ & $1.72\times 10^{-5}$ & $1.02\times10^{-5}$  \\   
\hline 
\end{tabular} 
\caption{Relative errors in $l^{\infty}$ norm at different times for the linear advection
problem. The  {\tt ode45} solver solution was used as a reference and
the solution was computed with a) the multirate method with cubic Hermite interpolator b) the  multirate method with linear interpolator  and c) the  single-rate TR-BDF2 method.}
\label{error_linearupwind}
\end{table}

%

The CPU times and total number of computed time steps are reported in Table \ref{SvsMhy} for
the multirate and single rate TR-BDF2 methods, respectively. For the multirate method, cubic interpolation was found to be necessary in order to achieve an improvement over the performance of the single rate method, while
linear interpolation lead to an excessive number of rejections of macro steps. 
In order to assess also the accuracy of the multirate approach, the relative errors obtained in this test case
are reported in Table \ref{error_linearupwind_ex}. More specifically, errors obtained by the multirate TR-BDF2  method with cubic Hermite interpolator, the  multirate TR-BDF2  method with linear interpolator  and the  single-rate TR-BDF2 method are compared. It can be seen that all three approaches yield comparable errors.
Since the first order upwind space discretization is responsible for the largest part of the error,
we also compute the error with respect to the reference solution obtained by discretizing in time  with the MATLAB {\tt  ode45}  solver set to a suitably stringent tolerance. The results, reported in
Table \ref{error_linearupwind}, show that also in this case the errors obtained with the three variants are
entirely analogous, albeit smaller by several orders of magnitude due the factoring out of the spatial discretization error.

\subsection{Burgers equation}
Finally, as a nonlinear test for hyperbolic equations, we consider the Riemann problem for the inviscid Burgers equation 

\begin{equation}
\dfrac{\partial u}{\partial t} + \dfrac{\partial}{\partial x} \left(\dfrac{1}{2}u^2\right) = 0, \qquad x\in[-1,3],\quad t>0,
\end{equation}
which amounts to assigning a piecewise constant  initial datum
$ u_0(x) = u_l \qquad x<0,  $
$ u_0(x) = u_r \qquad x>0. $
We first consider the case  $ u_l > u_r, $ for which a shock wave solution is obtained. 
The space discretization is given by   finite volume method with Rusanov numerical flux,
see e.g. \cite{leveque:2002}.
The interval time is $[0,1] $ with a number of cells equal to $400.$
The relative error tolerance  was
set to $ 10^{-4}, $ while the absolute error tolerance was set to
$ 10^{-6} $ and   the tolerance for the Newton solver required by the implicit TR-BDF2 method
was taken to be $ 10^{-8}$.
We took the initial size of the time step equal to $10^{-2}$. For the Riemann problem we used in this first case $u_l=1,$ while $u_r=0$.

\begin{figure}[!htb]
\centering
\advance\leftskip-1.cm   
\includegraphics[width=1.1\textwidth]{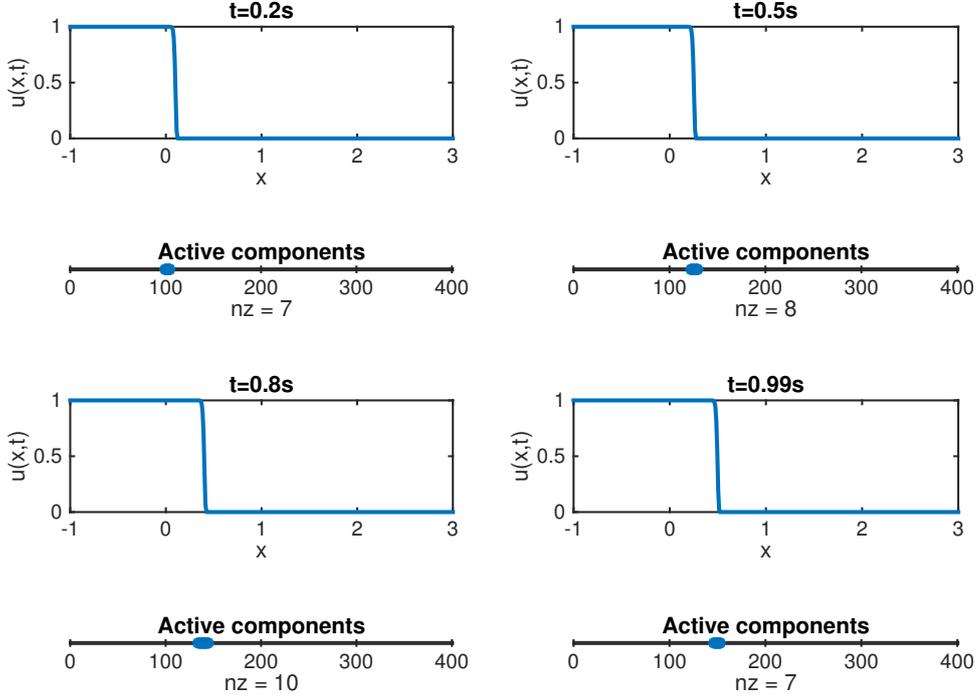}
\vspace*{-1.5cm}
\caption{Multirate TR-BDF2 integration with cubic interpolation for the Burgers equation Riemann problem with shock wave solution.}
\label{UBurgers1}
\end{figure}

\begin{figure}[!htb]
\centering   
\advance\leftskip-1.cm 
\includegraphics[width=0.7\textwidth]{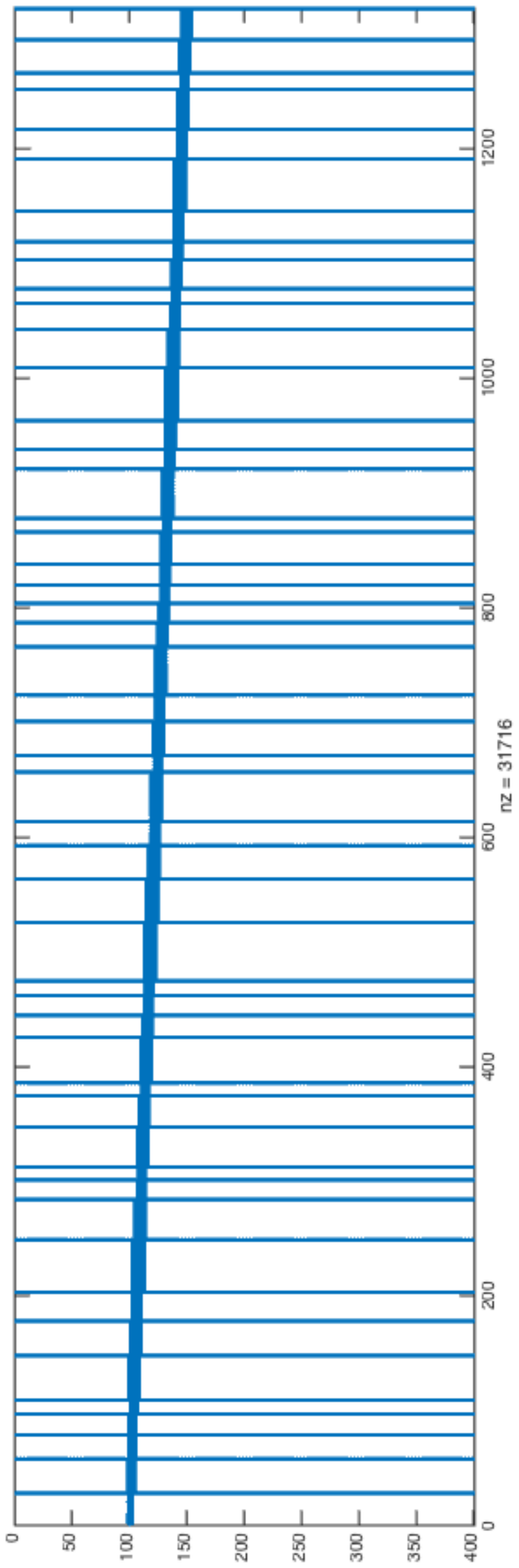}
\vspace*{-1.cm}
\caption{Space time diagram for the multirate TR-BDF2 algorithm for the Burgers equation Riemann problem with shock wave solution.}
\label{Burgers1mesh}
\end{figure}


\begin{table}[!htb]
\centering
\begin{tabular}{|c|c|c|}
\hline 
Method & CPU time  & Number of time steps \\ 
\hline 
Single rate TR-BDF2 & $78.23 s$ & $1008$ \\ 
\hline 
Multirate TR-BDF2  & $24.25 s$ & $1378$ \\ 
\hline  
\end{tabular} 
\caption{CPU time and number of time steps for the Burgers equation Riemann problem with shock wave solution.}
\label{Burgers1time}
\end{table}

\begin{table}[!htb]
\centering
\begin{tabular}{|c|p{2.5cm}| p{2.5cm}| p{2.5cm}| }
\hline
Time [s] & Cubic Interp. & Linear Interp. & Single rate \\ 
\hline
$0.2$ & $4.85\times10^{-1}$ & $4.83\times 10^{-1}$ &$4.81\times 10^{-1}$    \\ 
\hline
$0.5$ & $4.86\times 10^{-1}$&  $4.86\times 10^{-1}$ & $4.83\times 10^{-1}$    \\ 
\hline 
$0.8$ & $4.81\times10^{-1}$ & $4.82\times 10^{-1}$ & $4.80\times10^{-1}$   \\ 
\hline 
$0.99$ & $5.27\times10^{-1}$ & $5.31\times 10^{-1}$ & $5.24\times10^{-1}$   \\   
\hline
\end{tabular} 
\caption{Relative errors in $l^{\infty}$ norm at different times for the Burgers equation with shock wave solution.
The exact solution
was used as a reference  and
the solution was computed with a) the multirate method with cubic Hermite interpolator b) the  multirate method with linear interpolator  and c) the  single-rate TR-BDF2 method.}
\label{error_burgers1_ex}
\end{table}

\begin{table}[!htb]
\centering
\begin{tabular}{|c|p{2.5cm}| p{2.5cm}| p{2.5cm}|}
\hline
Time [s] & Cubic interp. & Linear interp. & Single rate   \\ 
\hline
$0.2$ & $4.37\times10^{-5}$ & $1.48\times 10^{-4}$ &$1.34\times 10^{-5}$    \\ 
\hline
$0.5$ & $2.76\times 10^{-4}$&  $6.08\times 10^{-4}$ & $3.5\times 10^{-5}$    \\ 
\hline 
$0.8$ & $6.65\times10^{-4}$ & $1.08\times 10^{-3}$ & $3.63\times10^{-5}$   \\ 
\hline 
$0.99$ & $9.18\times10^{-4}$ & $1.2\times 10^{-3}$ & $3.59\times10^{-5}$   \\   
\hline 
\end{tabular} 
\caption{Relative errors in $l^{\infty}$ norm at different times for the Burgers equation Riemann problem with shock wave solution.
The {\tt ode45} solution
was used as a reference  and
the solution was computed with a) the multirate method with cubic Hermite interpolator b) the  multirate method with linear interpolator  and c) the  single-rate TR-BDF2 method.}
\label{error_burgers1}
\end{table}

%
%

Figure \ref{UBurgers1} shows how the shock wave
solution evolves in space as time progresses.
It can be seen how the automatically selected set of
 active components moves along with the shock wave and does not increase in size.
 The same result is apparent from the space time diagram in Figure \ref{Burgers1mesh}.
 The CPU times and total number of computed time steps are reported in Table \ref{Burgers1time} for
the multirate and single rate TR-BDF2 methods, respectively. For the multirate method, cubic interpolation was found to be necessary in order to achieve an improvement over the performance of the single rate method, while
linear interpolation lead to an excessive number of rejections of macro steps. 
In order to assess also the accuracy of the multirate approach, the relative errors obtained in this test case
are reported in Table \ref{error_burgers1_ex}. More specifically, errors obtained by the multirate TR-BDF2  method with cubic Hermite interpolator, the  multirate TR-BDF2  method with linear interpolator  and the  single-rate TR-BDF2 method are compared. It can be seen that all three approaches yield comparable errors.
Since the first order upwind space discretization is responsible for the largest part of the error,
we also compute the error with respect to the reference solution obtained by discretizing in time  with the MATLAB {\tt  ode45}  solver set to a suitably stringent tolerance. The results, reported in
Table \ref{error_burgers1}, show that also in this case the errors obtained with the three variants are
entirely analogous, albeit smaller by several orders of magnitude due the factoring out of the spatial discretization error.

\begin{figure}[!htb]
\centering
\hspace*{-1.0cm}   
\includegraphics[width=0.45\textwidth]{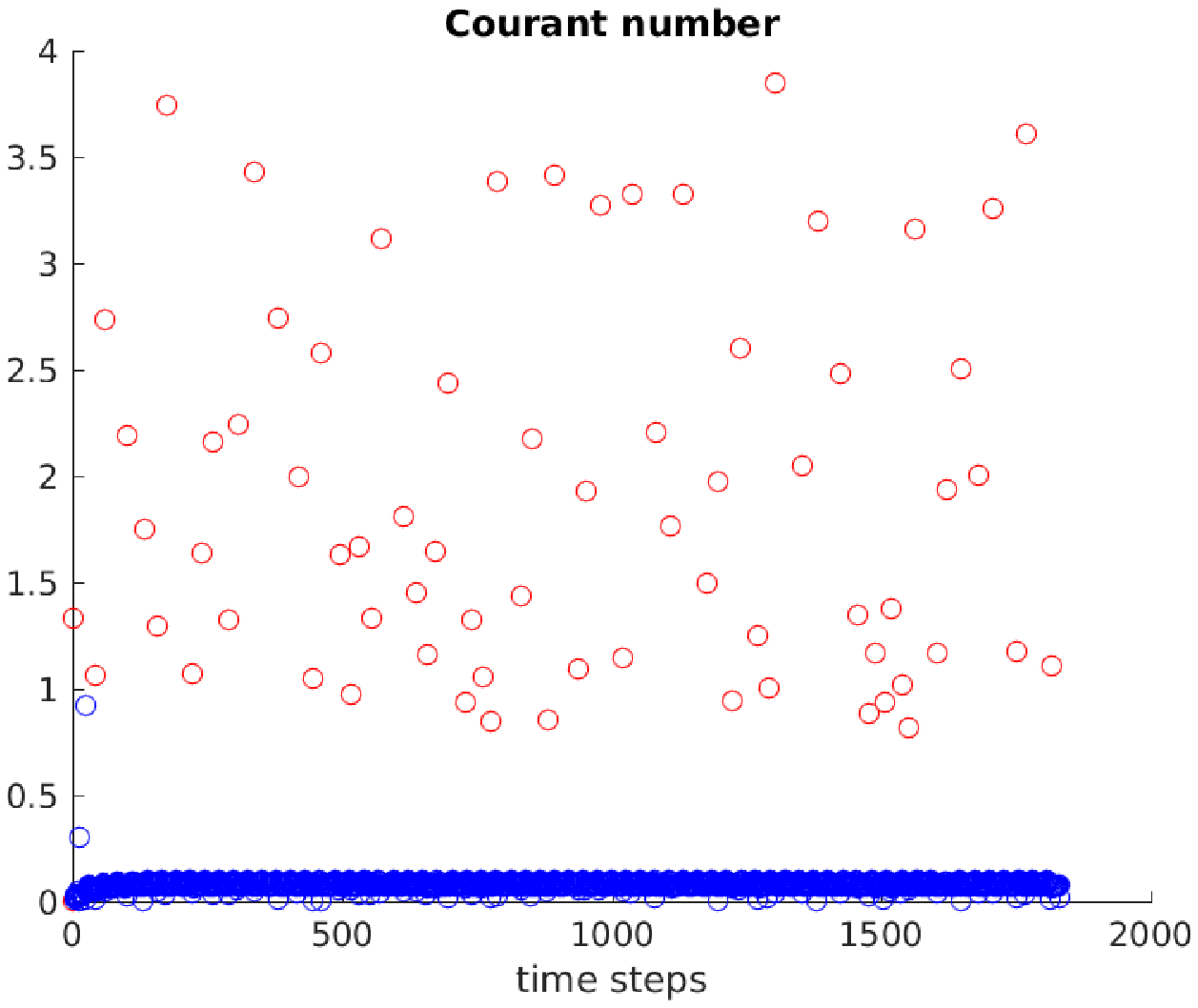} a)
\includegraphics[width=0.45\textwidth]{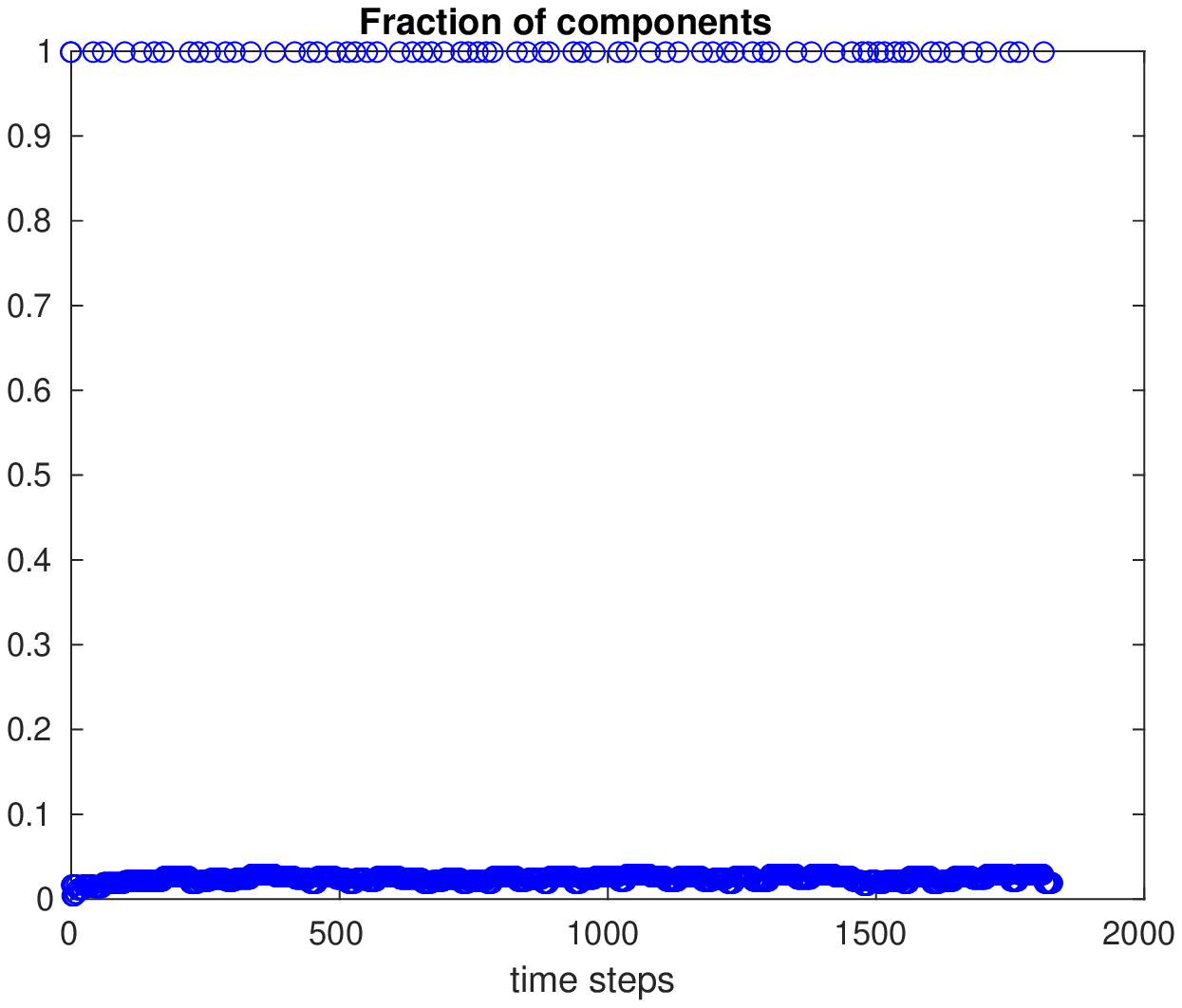}b)
\caption{a) Courant numbers for each time step and b) 
fraction of refined components for the Burgers equation with shock wave solution (red circles in a) denote  Courant numbers associated to the global time steps). }
\label{CourBurgers1}
\end{figure}

In order to highlight the capability of the adaptive strategy to select automatically a time step that
is compatible with the underlying dynamics, we 
  have also reported the maximum Courant numbers at each time step, computed as
\begin{equation}
\max_{i=1,...,N_{x}}|f'(c_{i}^{n})|\dfrac{h_n}{\Delta x}   
\end{equation}
where $N_x$ is the total number of cells. In Figure \ref{CourBurgers1} a),
the Courant number values associated to the global time steps are represented by red circles,
while those associated to local, refined time steps are represented by blue symbols.
It can be seen that, while   large Courant numbers are used for the latent components,
much smaller values are employed for the active components that correspond to the locations
actually crossed by the shock wave.  In Figure \ref{CourBurgers1} b) the corresponding fractions of active
components are displayed, highlighting the significant reduction in computational cost in this case.


 We then consider the case $u_l > u_r, $ for which a rarefaction wave solution is obtained.
More specifically, we consider $u_l=0$ and $u_r=1$. The other parameters are the same as in the previous case.

\begin{figure}[!htb]
\centering
\hspace*{-1.cm}   
\includegraphics[width=1.1\textwidth]{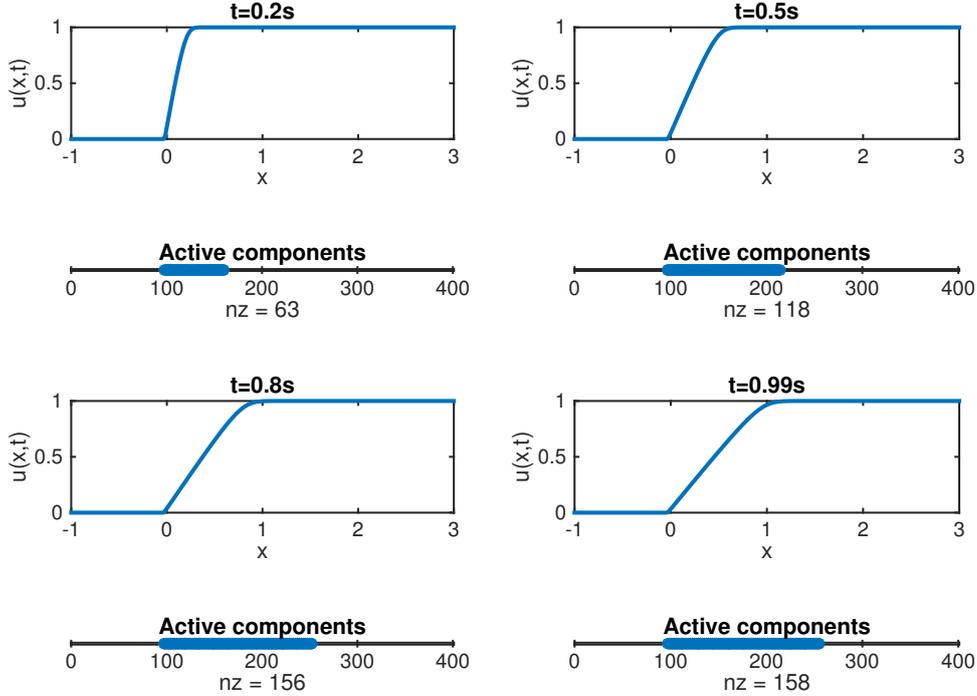}
\vspace*{-2.cm}
\caption{Multirate TR-BDF2 integration with cubic interpolation for the Burgers equation Riemann problem with 
 rarefaction wave solution.}
\label{UBurgers2}
\end{figure}

\begin{figure}[!htb]
\centering   
\advance\leftskip-1.cm 
\includegraphics[width=0.8\textwidth]{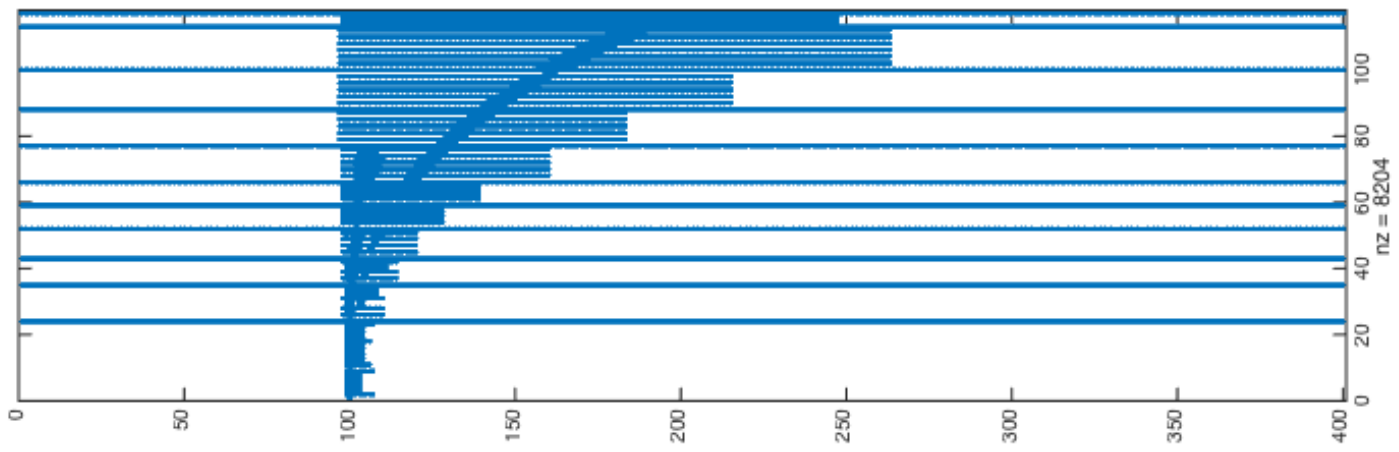}
\vspace*{-1.cm}
\caption{Space time diagram for the multirate TR-BDF2 algorithm for the Burgers equation Riemann problem with rarefaction wave solution.}
\label{Burgers2mesh}
\end{figure}


\begin{table}[!htb]
\centering
\begin{tabular}{|c|c|c|}
\hline 
Method & CPU time  & Numb. time steps \\ 
\hline 
Single rate TR-BDF2 & $7.13 s$ & $79$ \\ 
\hline 
Multirate TR-BDF2 & $2.99 s$ & $169$ \\ 
\hline 
\end{tabular} 
\caption{CPU time and number of time steps for the rarefaction wave.}
\label{Burgers2time}
\end{table}

\begin{table}[!htb]
\centering 
\begin{tabular}{|c|p{2.5cm}| p{2.5cm}| p{2.5cm}| }
\hline
Time [s] & Cubic interp.  & Linear interp. & Single rate \\ 
\hline
$0.2$ & $3.59\times10^{-1}$ & $3.59\times 10^{-1}$ &$3.59\times 10^{-1}$ \\ 
\hline
$0.5$ & $3.16\times 10^{-1}$&  $3.16\times 10^{-1}$ & $3.14\times 10^{-1}$ \\ 
\hline 
$0.8$ & $2.93\times10^{-1}$ & $2.95\times 10^{-1}$ & $2.85\times10^{-1}$\\ 
\hline 
$0.99$ & $2.84\times10^{-1}$ & $2.92\times 10^{-1}$ & $2.83\times10^{-1}$ \\   
\hline 
\end{tabular} 
\caption{Relative errors in $l^{\infty}$ norm at different times for the Burgers equation Riemann problem with rarefaction wave solution.
The exact solution
was used as a reference  and
the solution was computed with a) the multirate method with cubic Hermite interpolator b) the  multirate method with linear interpolator  and c) the  single-rate TR-BDF2 method. }
\label{error_burgers2_ex}
\end{table}
%
%

\begin{table}[!htb]
\centering
\begin{tabular}{|c|p{2.5cm}| p{2.5cm}| p{2.5cm}|  }
\hline
Time [s] & Cubic interp.  & Linear interp. & Single rate  \\ 
\hline
$0.2$ & $3.13\times10^{-4}$ & $4.42\times 10^{-4}$ &$1.92\times 10^{-4}$  \\ 
\hline
$0.5$ & $5.53\times 10^{-4}$&  $9.91\times 10^{-4}$ & $6.99\times 10^{-4}$  \\ 
\hline 
$0.8$ & $7.57\times10^{-4}$ & $1.75\times 10^{-3}$ & $7.25\times10^{-4}$ \\ 
\hline 
$0.99$ & $7.48\times10^{-4}$ & $1.25\times 10^{-3}$ & $6.64\times10^{-4}$  \\   
\hline 
\end{tabular} 
\caption{Relative errors in $l^{\infty}$ norm at different times for the Burgers equation Riemann problem with rarefaction wave solution.
The {\tt ode45} solution
was used as a reference  and
the solution was computed with a) the multirate method with cubic Hermite interpolator b) the  multirate method with linear interpolator  and c) the  single-rate TR-BDF2 method.}
\label{error_burgers2}
\end{table}

Figure \ref{UBurgers2} shows how the rarefaction wave
solution evolves in space as time progresses.
It can be seen how the automatically selected set of
 active components moves along with the   wave.
 The same result is apparent from the space time diagram in Figure \ref{Burgers2mesh}.
 Since in this case the solution  undergoes significant changes in a  larger number of cells,
 the number of active components is correspondingly larger than in the shock wave case.

 The CPU times and total number of computed time steps are reported in Table \ref{Burgers2time} for
the multirate and single rate TR-BDF2 methods, respectively. For the multirate method, cubic interpolation was found to be necessary in order to achieve an improvement over the performance of the single rate method, while
linear interpolation lead to an excessive number of rejections of macro steps. 
In order to assess also the accuracy of the multirate approach, the relative errors obtained in this test case
are reported in Table \ref{error_burgers2_ex}. More specifically, errors obtained by the multirate TR-BDF2  method with cubic Hermite interpolator, the  multirate TR-BDF2  method with linear interpolator  and the  single-rate TR-BDF2 method are compared. It can be seen that all three approaches yield comparable errors.
Since the first order upwind space discretization is responsible for the largest part of the error,
we also compute the error with respect to the reference solution obtained by discretizing in time  with the MATLAB {\tt  ode45}  solver set to a suitably stringent tolerance. The results, reported in
Table \ref{error_burgers2}, show that also in this case the errors obtained with the three variants are
entirely analogous, albeit smaller by several orders of magnitude due the factoring out of the spatial discretization error.

\begin{figure}[!htb]
\centering
\hspace*{-1.0cm}   
\includegraphics[width=0.45\textwidth]{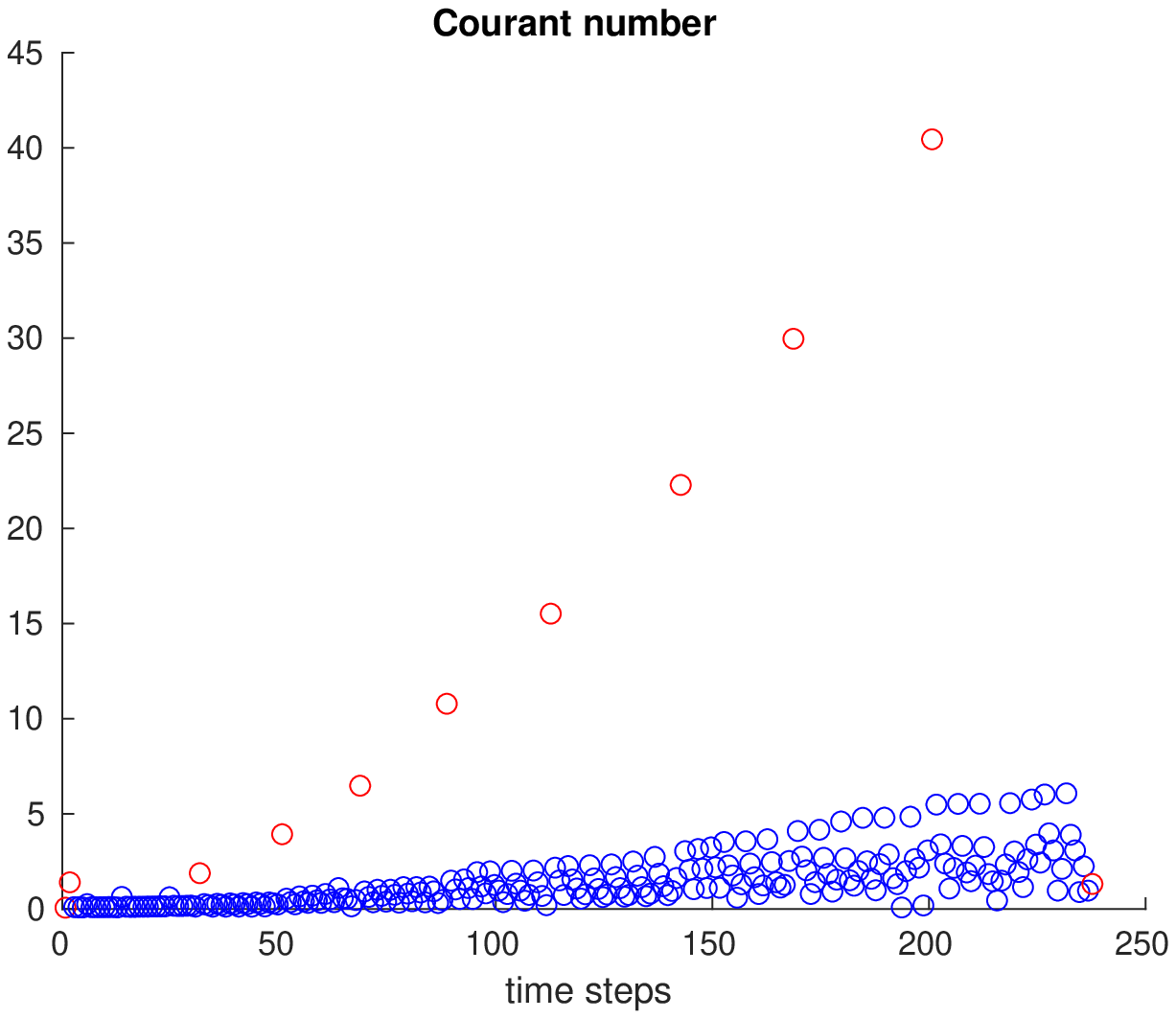} a)
\includegraphics[width=0.45\textwidth]{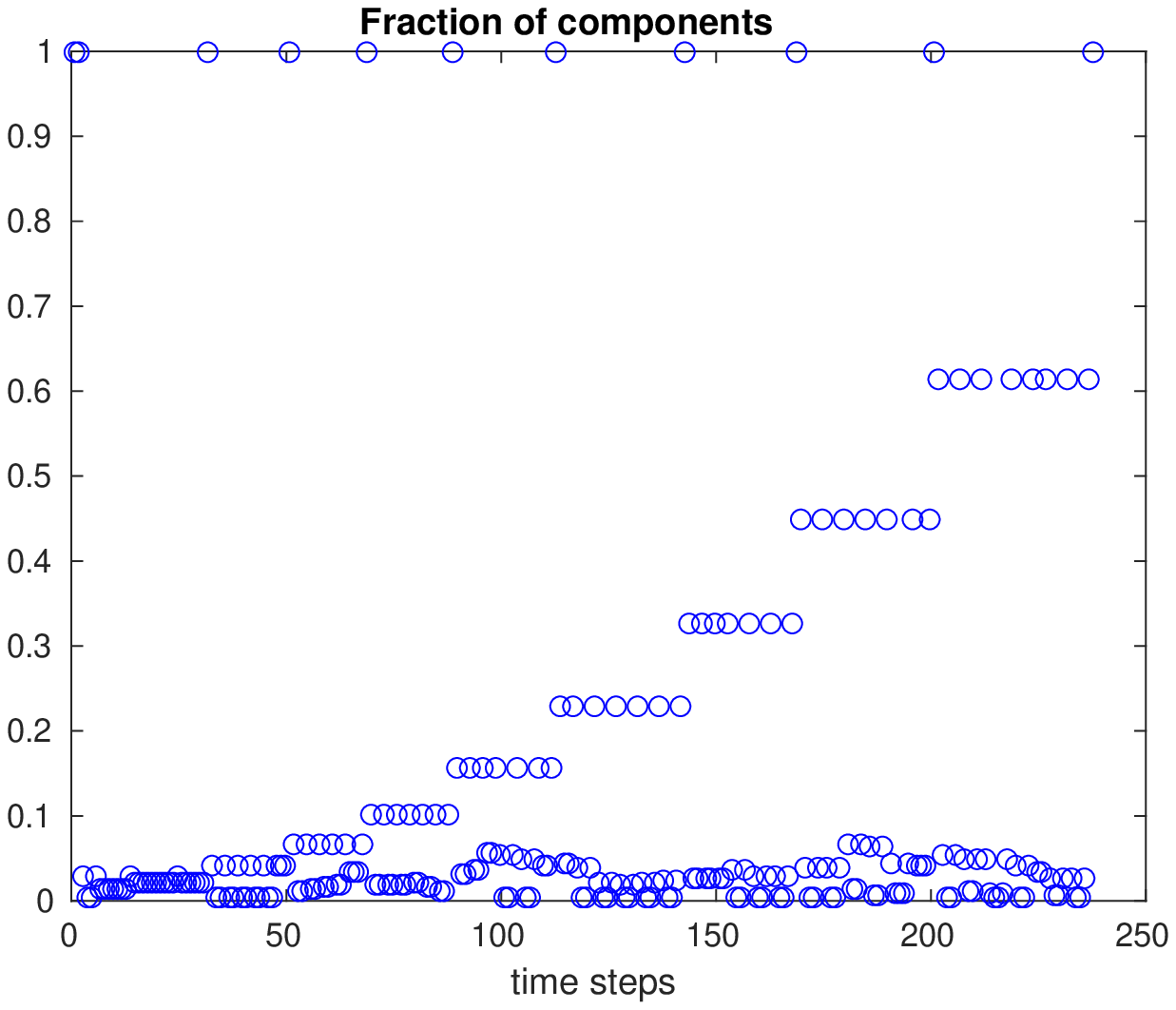}b)
\caption{a) Courant numbers for each time step and b) 
fraction of refined components for the Burgers equation with rarefaction wave solution (red circles in a) denote  Courant numbers associated to the global time steps). }
\label{CourBurgers2}
\end{figure}

In Figure \ref{CourBurgers2} a)  we have reported the Courant number values associated to each time step. Comparing to the shock wave case, we can see that larger Courant numbers
 are allowed by the error estimator also in the case of refined time steps, as a consequence of the different dynamics of this kind of solution.  In Figure \ref{CourBurgers2} b) the corresponding fractions of active
components are displayed, highlighting the lesser reduction in computational cost in this case with respect
to the shock wave solution.

%

 \section{Conclusions}
 \label{conclu} \indent

 We have  proposed a self-adjusting multirate method  
that relies on the TR-BDF2 method as fundamental single rate solver. 
 The TR-BDF2 method has  a number  of interesting properties, such as cheap error estimation
 via an embedded third order method and  continuous output via a cubic Hermite interpolant, that can
 be exploited to increase efficiency and accuracy of a multirate approach. 
 
 Our  self-adjusting multirate TR-BDF2 method  has been coupled to a
    partitioning and time step selection criterion based on the technique proposed in  \cite{fok:2015}.  
  The stability of the method has been analyzed in the framework of the classical linear model problem.
  Even though the results achieved were only based on the numerical computation of the stability function norm,
  we show that  for a range of relevant  stiff model problems  the method has remarkable stability properties.
  This is true not only for contractive problems with strictly dissipative behaviour, but also for problems
  with purely  imaginary eigenvalues, which suggests that the method could be advantageous also for
  the application to hyperbolic PDEs and structural mechanics problems.
  The numerical results obtained on several standard benchmarks show that application of the proposed method leads  to significant efficiency gains, that are analogous to those achieved by other self-adjusting approaches when compared   to their corresponding single rate variants. In particular,  the has been applied to the time discretization of nonlinear, hyperbolic partial differential equations, allowing to achieve automatic detection
  of complex localized phenomena such as shock waves and significant reductions in computational cost.
  
  In a companion work, we will focus on the extensive application of the method described in this paper
  to nonlinear hyperbolic equations and on the derivation of mass conservative versions
  of  this multirate approach.

\section*{Acknowledgements} This paper contains developments of the results 
presented in the PhD thesis by A.R. \cite{ranade:2016}  and in the Master thesis by L.D. \cite{delpopolo:2015}, both at Politecnico di Milano.
 A.R.'s PhD grant was supported by the Erasmus Mundus Heritage project, while he is presently supported by the Science Foundation Ireland grant 13/RC/2094.
L.B. was partially supported by the INDAM - GNCS 2015 project {\it Metodi numerici semi-impliciti e semi-Lagrangiani  per sistemi iperbolici di leggi di bilancio.}
 Useful discussions with Luca Formaggia are gratefully acknowledged, as well the comments by
  Nicola Guglielmi and Claus F\"uhrer on the first draft of the thesis by Akshay Ranade.

\bibliographystyle{plain}
\bibliography{multirate_new}

\end{document}